\documentclass[10pt]{amsart}

\usepackage{amsmath,amssymb,amsfonts,epsfig}

\newcommand\risS[6]{\raisebox{#1pt}[#5pt][#6pt]{\begin{picture}(#4,15)(0,0)
  \put(0,0){\includegraphics[width=#4pt]{#2.eps}} #3
     \end{picture}}}
\newcommand\smf[1]{\risS{-4}{#1}{}{15}{15}{8}}
\newcommand\laf[1]{\begin{picture}(65,50)(0,0)
          \put(0,0){\includegraphics[width=65pt]{#1.eps}}
                   \end{picture}}
\def\a{\alpha}
\def\b{\beta}
\def\d{\delta}
\def\ve{\varepsilon}
\def\Ga{\Gamma}
\def\cS{\mathcal S}
\def\cF{\mathcal F}
\def\bc{\mathrm{bc}}
\newtheorem{thm}{Theorem}[section]
\newtheorem{defn}[thm]{Definition}
\newtheorem{Example}[thm]{Example}
\def\kb#1{[ #1 ]}
\def\wo{\overline}

\begin{document}

\title[Thistlethwaite's theorem for virtual links]
    {Thistlethwaite's theorem for virtual links}
\author[Sergei~Chmutov]{Sergei~Chmutov}
\author[Jeremy Voltz]{Jeremy Voltz}
\date{}


\begin{abstract}
The celebrated Thistlethwaite theorem relates the Jones polynomial of a
link with the Tutte polynomial of the corresponding planar graph.
We give a generalization of this theorem to virtual links. In this case,
the graph will be embedded into a (higher genus) surface. For such graphs
we use the generalization of the Tutte polynomial discovered by
B.~Bollob\'as and O.~Riordan.
\end{abstract}

\maketitle

\section*{Introduction} \label{s:intro}
Regions of a link diagram can be colored black and white in a checkerboard pattern. 
Putting a vertex in each black region and connecting two vertices by an edge if the 
corresponding regions share a crossing yields a planar graph. In 1987 Thistlethwaite \cite{Th} proved that the Jones polynomial of an alternating link can be obtained as a specialization of the Tutte polynomial of the corresponding planar graph. L.~Kauffman \cite{Ka2} generalized the theorem to arbitrary links using signed graphs and extending the Tutte polynomial to them. 
An expression for the Jones polynomial in terms of the Bollob\'as-Riordan polynomial,
without signed graphs, was found in \cite{DFKLS}. The idea to use 
the Bollob\'as-Riordan polynomial instead of the Tutte polynomial belongs to Igor Pak. It was first realized in \cite{CP}, where Thistlethwaite's theorem was generalized to 
checkerboard colorable virtual links.

Here we shall generalize this theorem to {\it arbitrary} virtual links. 

We recall the basic definitions of virtual links, their Jones polynomial through the Kauffman bracket, ribbon graphs, and the Bollob\'as-Riordan polynomial
in Sections \ref{s:vl} and \ref{s:rg}. The key construction of a ribbon graph from
a digram of a virtual link is explained in Section \ref{s:vl2rg}. Our main theorem is formulated and proved in Section \ref{s:mt}.

The work has been done as part of the Summer 2006 VIGRE working group ``Knots and Graphs" 
(\verb#http://www.math.ohio-state.edu/~chmutov/wor-gr-su06/wor-gr.htm#)
at the Ohio State University, funded by NSF grant DMS-0135308. We are grateful to 
O.~Dasbach and N.~Stoltzfus for useful and stimulating conversations.

\section{Virtual links and the Kauffman bracket}\label{s:vl}

The theory of \emph{virtual links} was discovered independently by
L.~Kauffman~\cite{Ka3} and M.~ Goussarov, M.~Polyak, and O.~Viro~\cite{GPV} around 1998.
According to Kauffman's approach, virtual links are represented by diagrams similar to ordinary knot diagrams, except some crossings are designated as {\it virtual}. Virtual crossings should be understood not as
crossings but rather as defects of our two-dimensional pictures. They should
be treated in the same way as the extra crossings appearing in planar pictures of
non-planar graphs. Here are some examples.
$$\risS{-18}{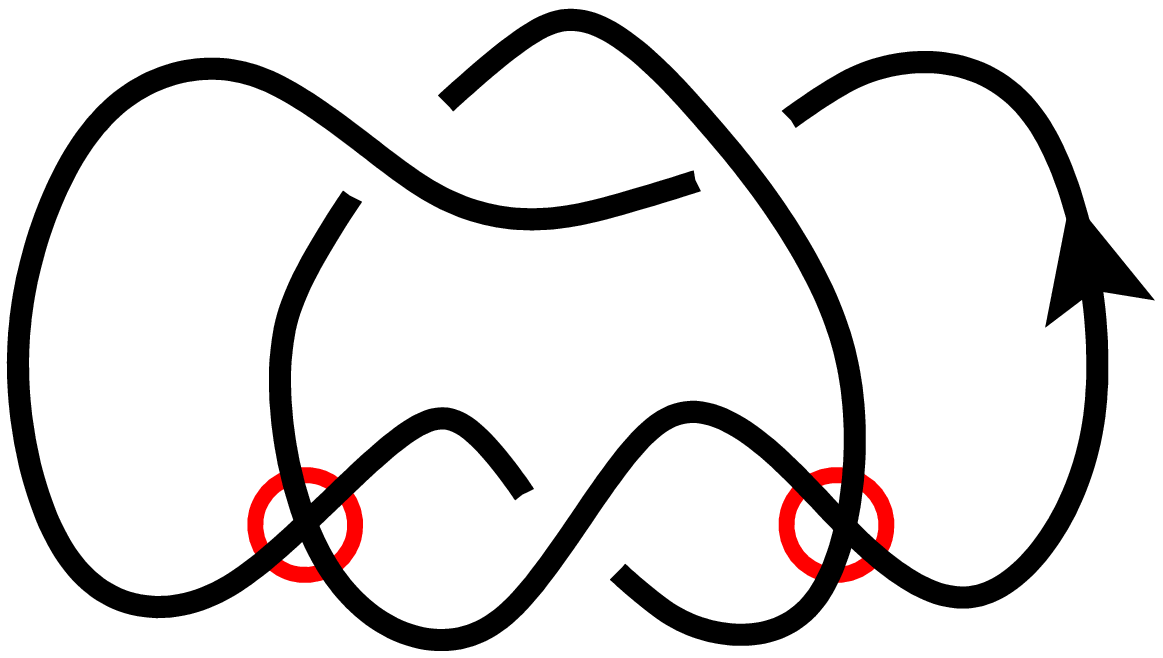}{}{65}{20}{20}\hspace{3cm}
  \risS{-18}{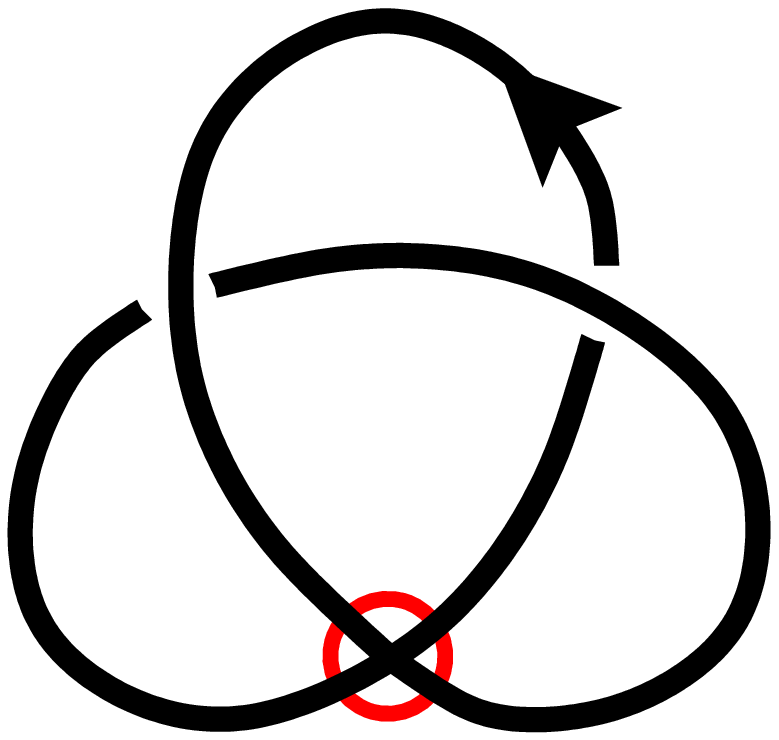}{}{45}{0}{20}\hspace{3cm}
  \risS{-18}{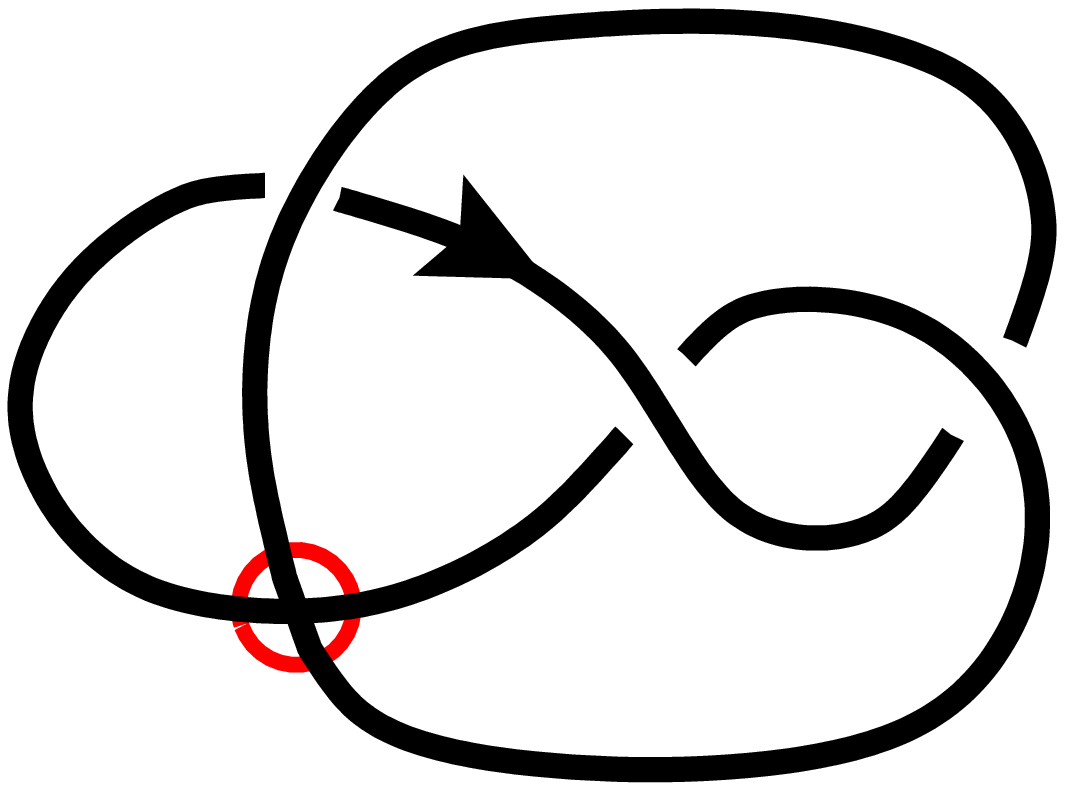}{}{55}{0}{20}
$$
The virtual crossings in these pictures are circled to distinguish them from the classical ones.

Virtual link diagrams are considered up to the {\it classical}
Reidemeister moves involving classical crossings:
$$\risS{-18}{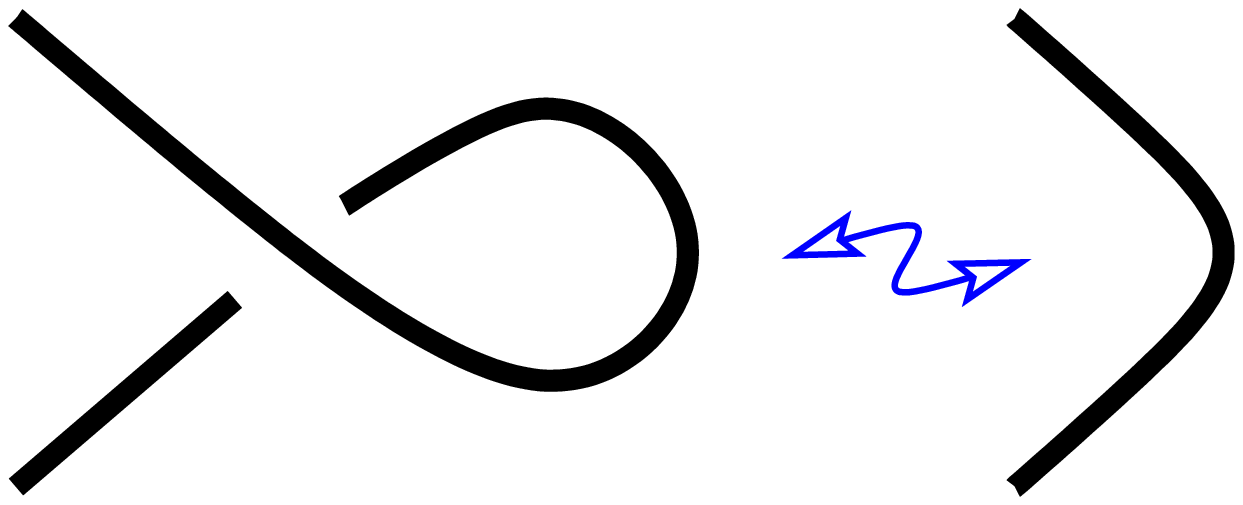}{}{75}{20}{20}\qquad\qquad
  \risS{-18}{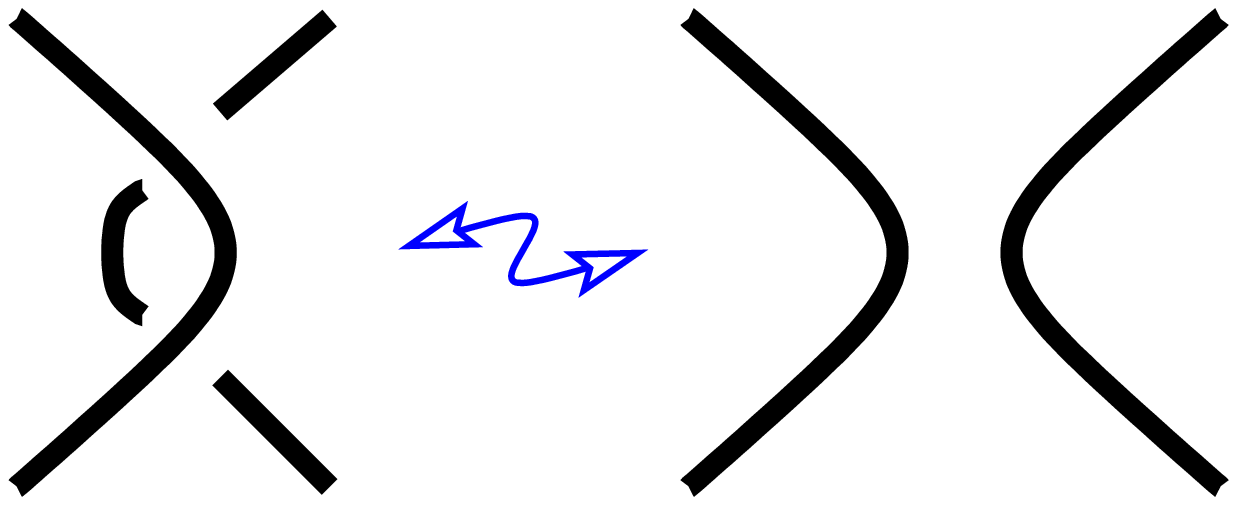}{}{75}{0}{20}\qquad\qquad
  \risS{-18}{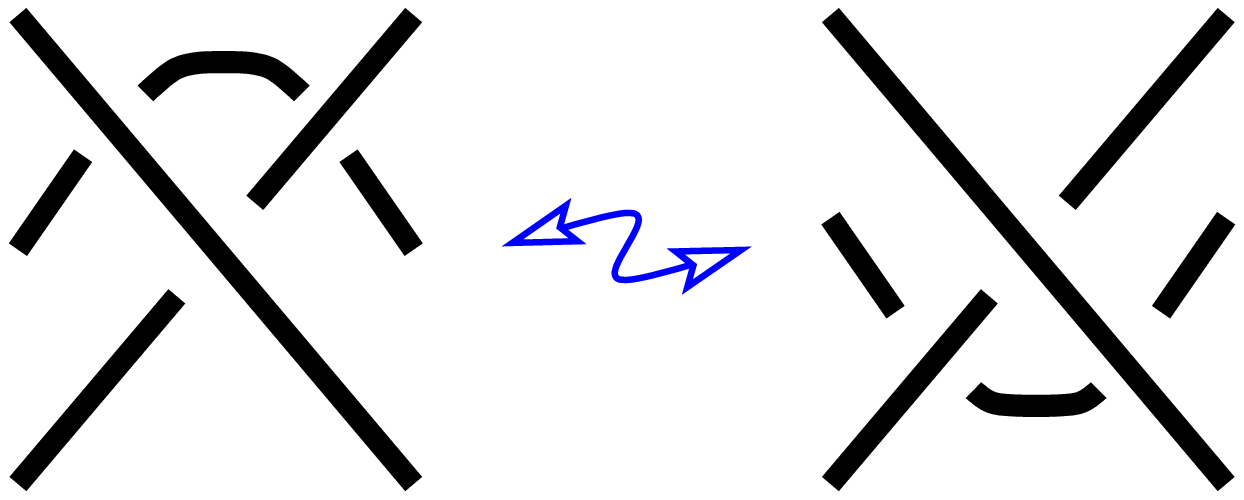}{}{75}{0}{20}
$$
and the {\it virtual} Reidemeister moves:
$$\risS{-18}{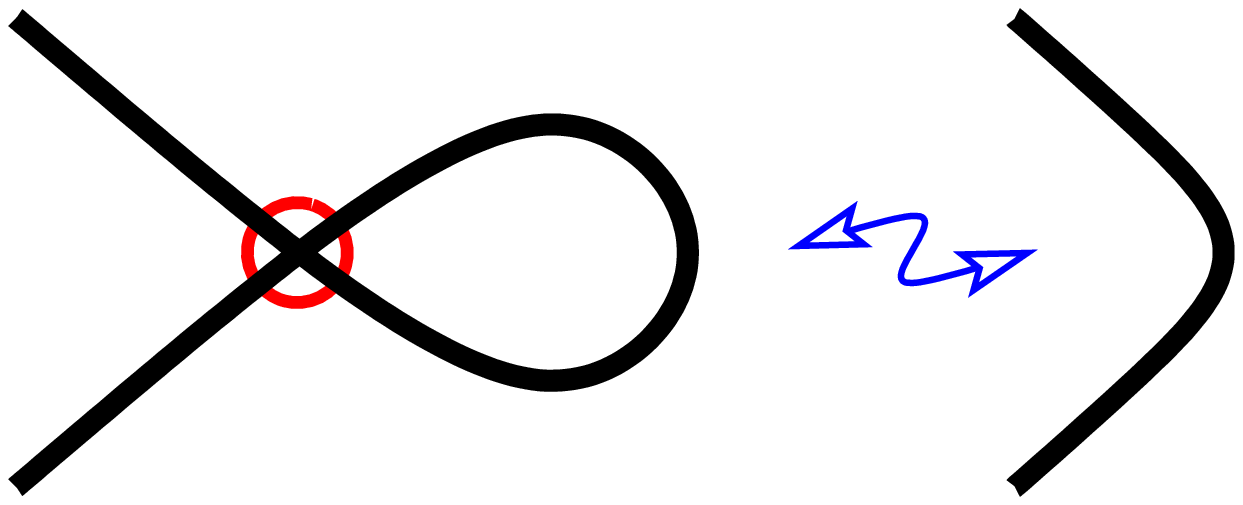}{}{75}{20}{20}\qquad\quad
  \risS{-18}{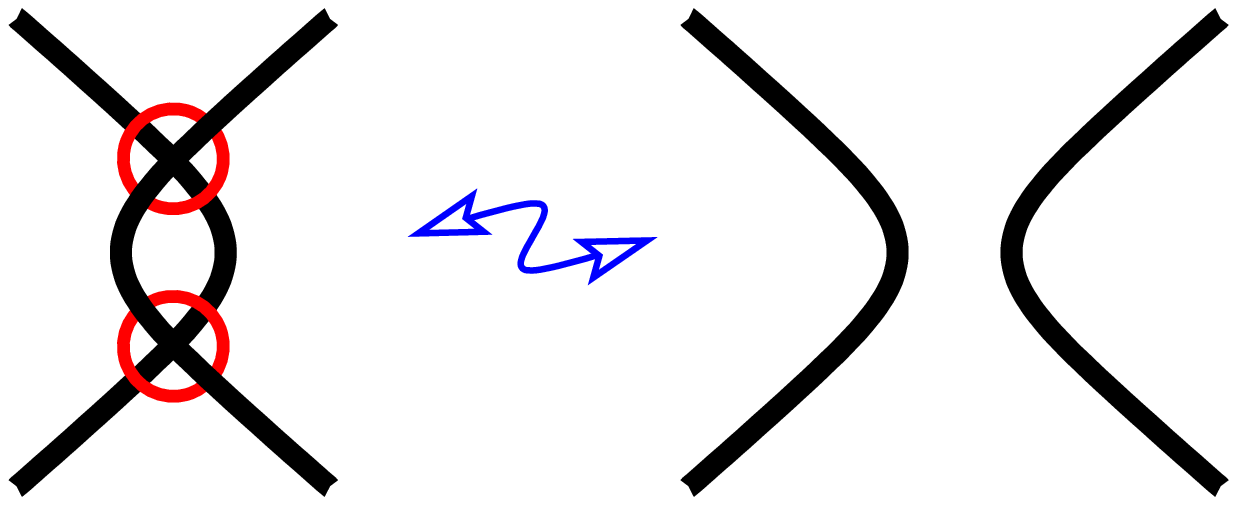}{}{75}{0}{20}\qquad\quad
  \risS{-18}{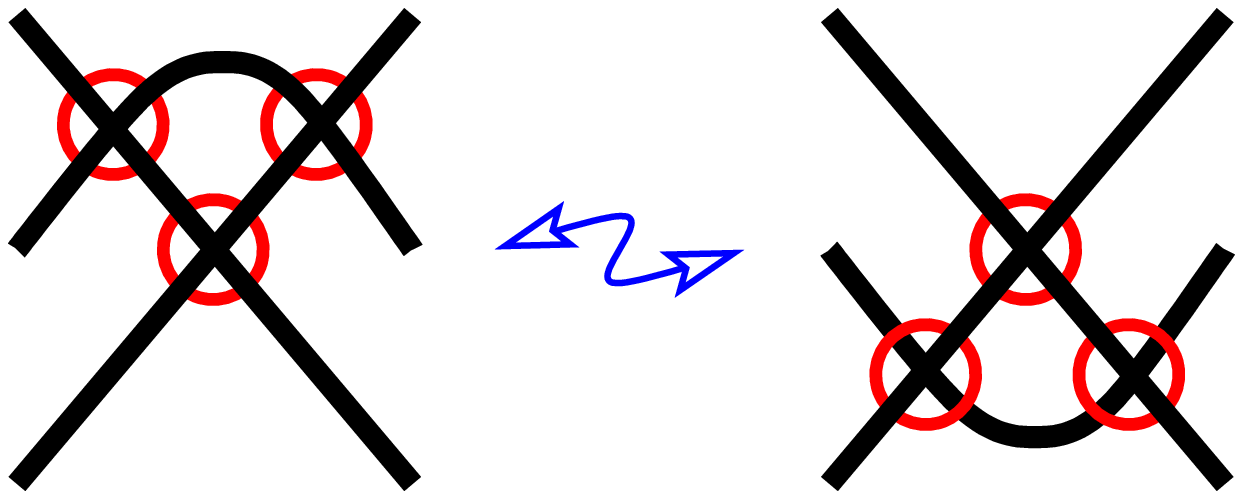}{}{75}{0}{20}\qquad\quad
  \risS{-18}{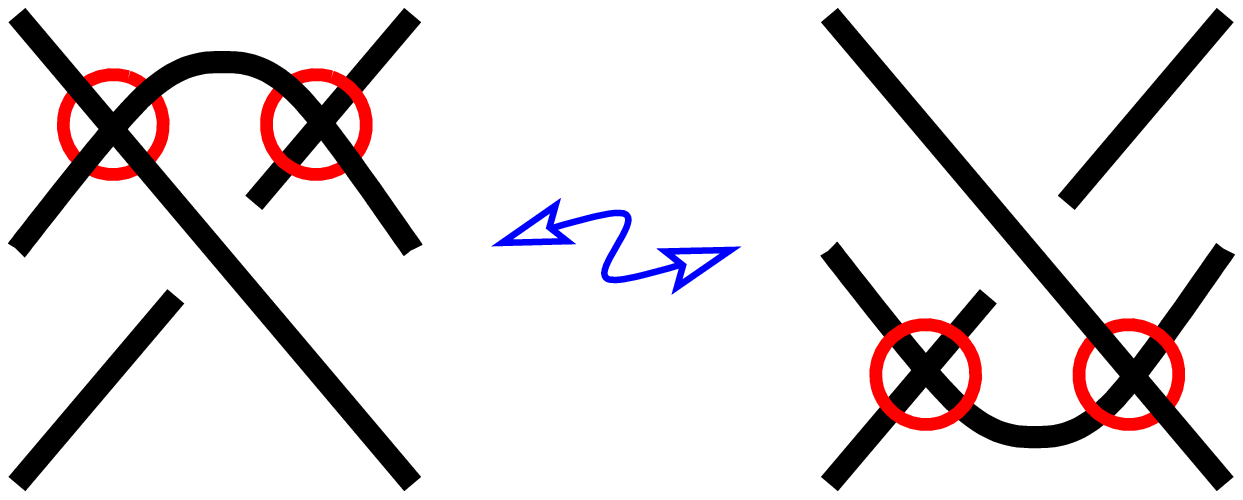}{}{75}{0}{30}
$$

The Kauffman bracket and the Jones polynomial for virtual links are defined in the same way as for classical ones. Let $L$ be a virtual link diagram.
Consider two ways of resolving a classical crossing.
The {\it $A$-splitting},\ $\smf{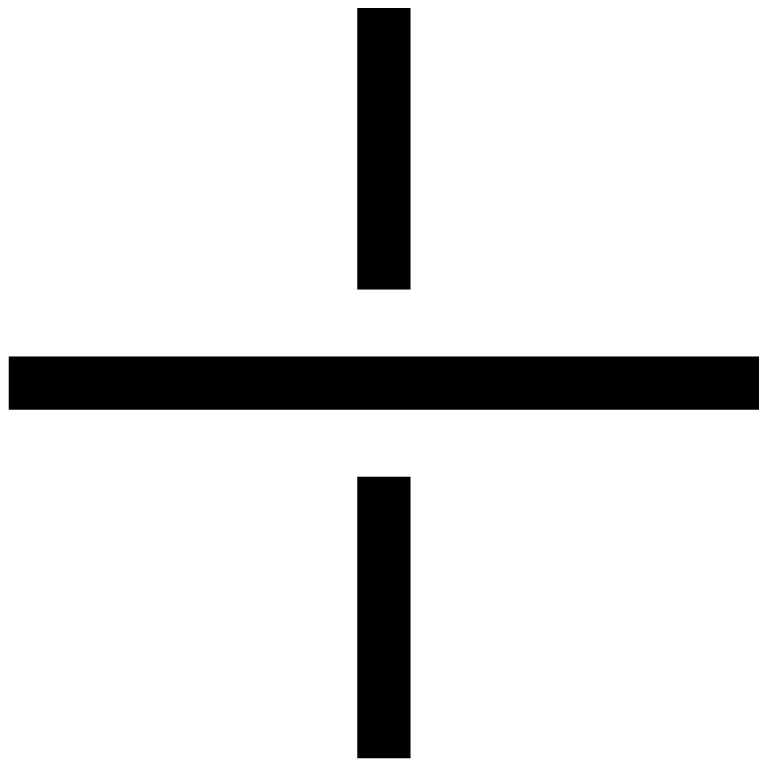}\ \leadsto\ \smf{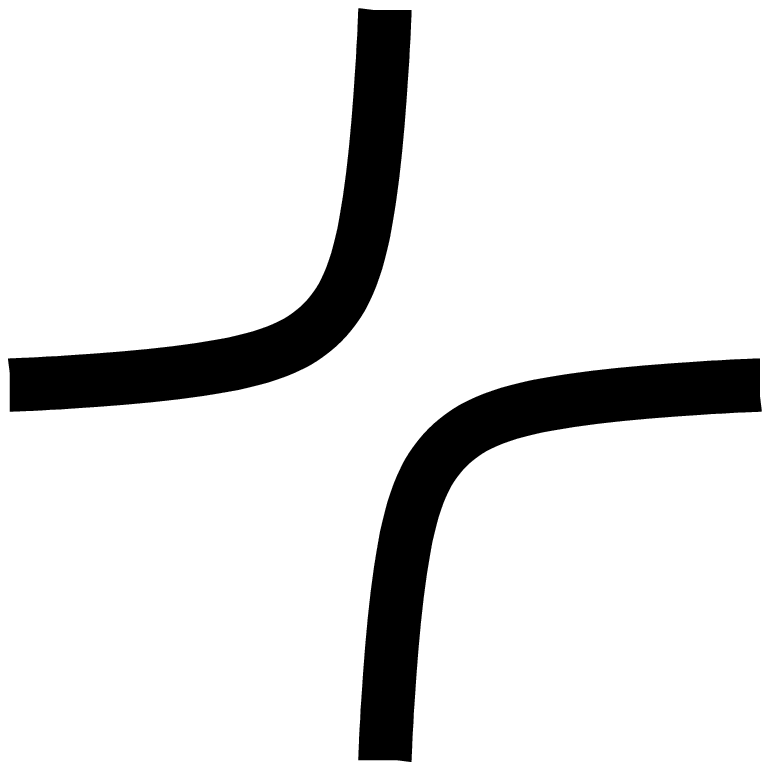}$\ ,
is obtained by joining the two vertical angles swept out by the overcrossing arc when
it is rotated counterclockwise toward the undercrossing arc.
Similarly, the {\it $B$-splitting},\ $\smf{cr}\ \leadsto\ \smf{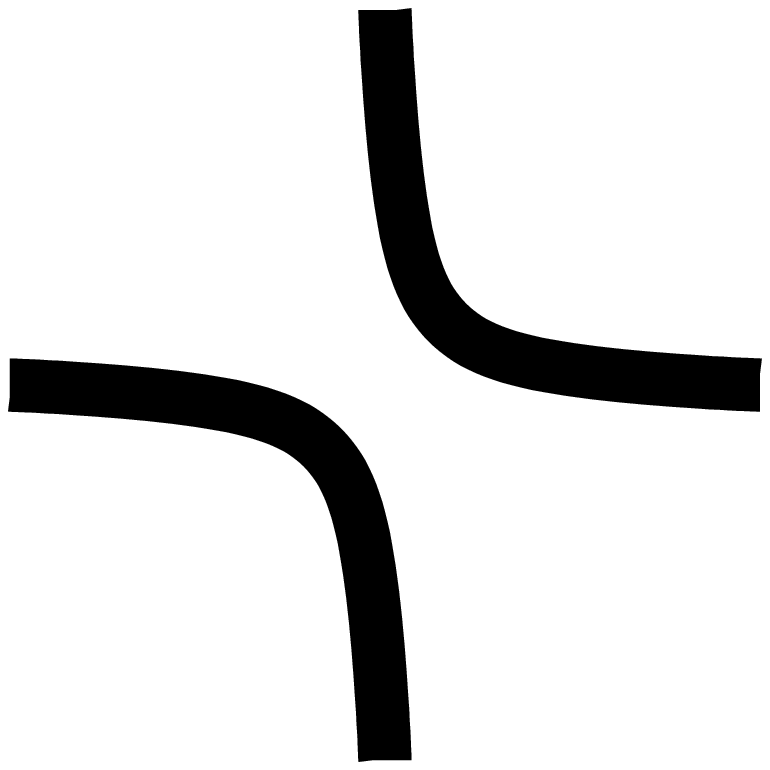}$\ , is
obtained by joining the other two vertical angles. A {\it state} $S$ of
a link diagram $L$
is a choice of either an $A$- or $B$-splitting at each classical crossing of the diagram.
Denote by $\cS(L)$ the set of the states of $L$.
Clearly, a diagram $L$ with $n$ crossings has $|\cS(L)| = 2^n$
different states.

Denote by $\a(S)$ and $\b(S)$ the numbers of $A$-splittings and $B$-splittings
in a state $S$, respectively.  Also, denote by $\d(S)$ the number of
components of the curve obtained from the link
diagram $L$ by 
splitting according to the state $S \in \cS(L)$. Note that virtual crossings do not connect components.

\begin{defn}\label{def:kb}
The \emph{Kauffman bracket} of a diagram $L$ is a polynomial in three variables
$A$, $B$, $d$ defined by the formula:
$$ \kb{L} (A,B,d)\ :=\ \sum_{S \in \cS(L)} \,
A^{\a(S)} \, B^{\b(S)} \, d^{\,\d(S)-1}\,.
$$
\end{defn}

Note that $\kb{L}$ is \emph{not} a topological
invariant of the link and in fact depends on the link
diagram. However, it defines the \emph{Jones polynomial}
$J_L(t)$ by a simple substitution:
$$J_L(t)\, := (-1)^{w(L)} t^{3w(L)/4} \kb{L} (t^{-1/4}, t^{1/4}, -t^{1/2}-t^{-1/2})\ .
$$
Here $w(L)$ denotes the {\it writhe}, determined
by orienting $L$ and taking the sum over the classical
crossings of $L$ of the following signs\,:
$$\begin{picture}(50,30)(0,0)
  \put(0,0){\includegraphics[width=50pt]{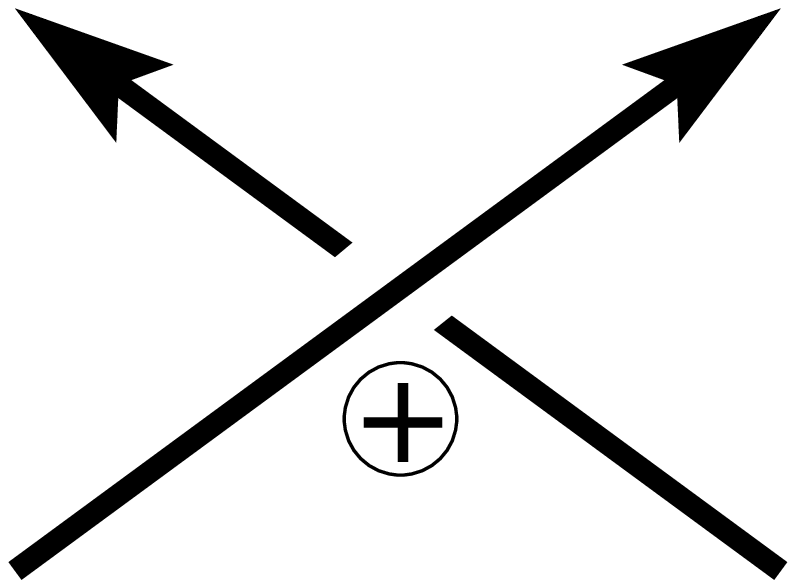}}
     \end{picture} \hspace{3cm}
\begin{picture}(50,40)(0,0)
  \put(0,0){\includegraphics[width=50pt]{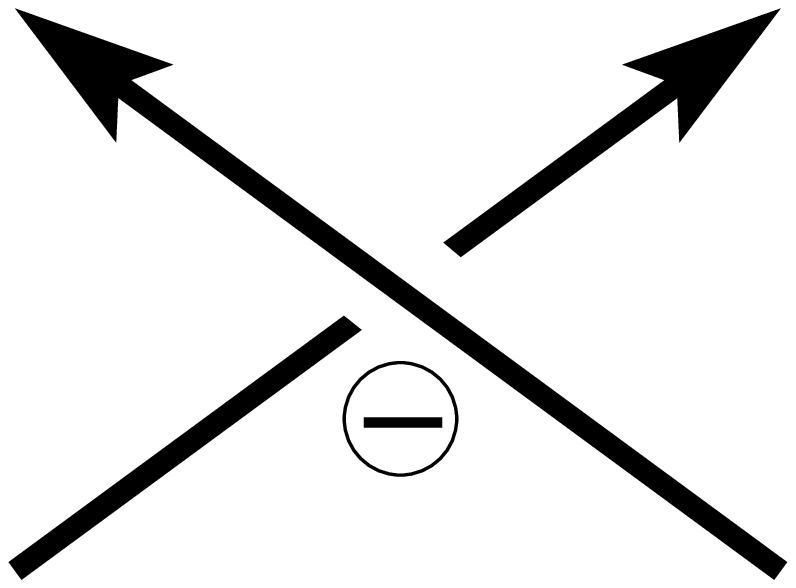}}
     \end{picture}
$$
The Jones polynomial is a classical topological invariant (see e.g.~\cite{Ka1}).

\begin{Example}\label{ex1} {\rm
Consider the third virtual knot diagram
$L$ from the example above. It is shown on the left of the table below.
It has one virtual and three classical crossings (one positive and two negative). 
So there are eight states and $w(L)=-1$.
The curves obtained by the splittings and the corresponding
parameters $\a(S)$, $\b(S)$, and $\d(S)$ are shown in the
remaining columns of the table.
$$\label{kb-table}
\begin{array}{c||c|c|c|c}
\laf{v41} & \laf{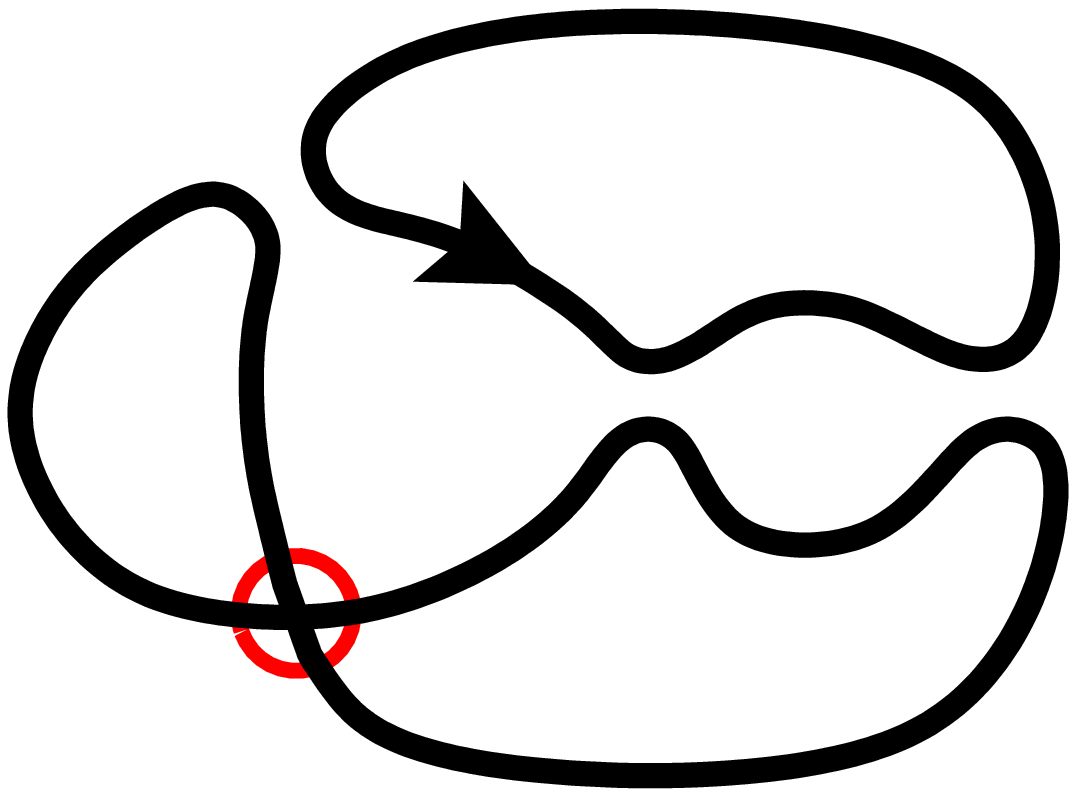} & \laf{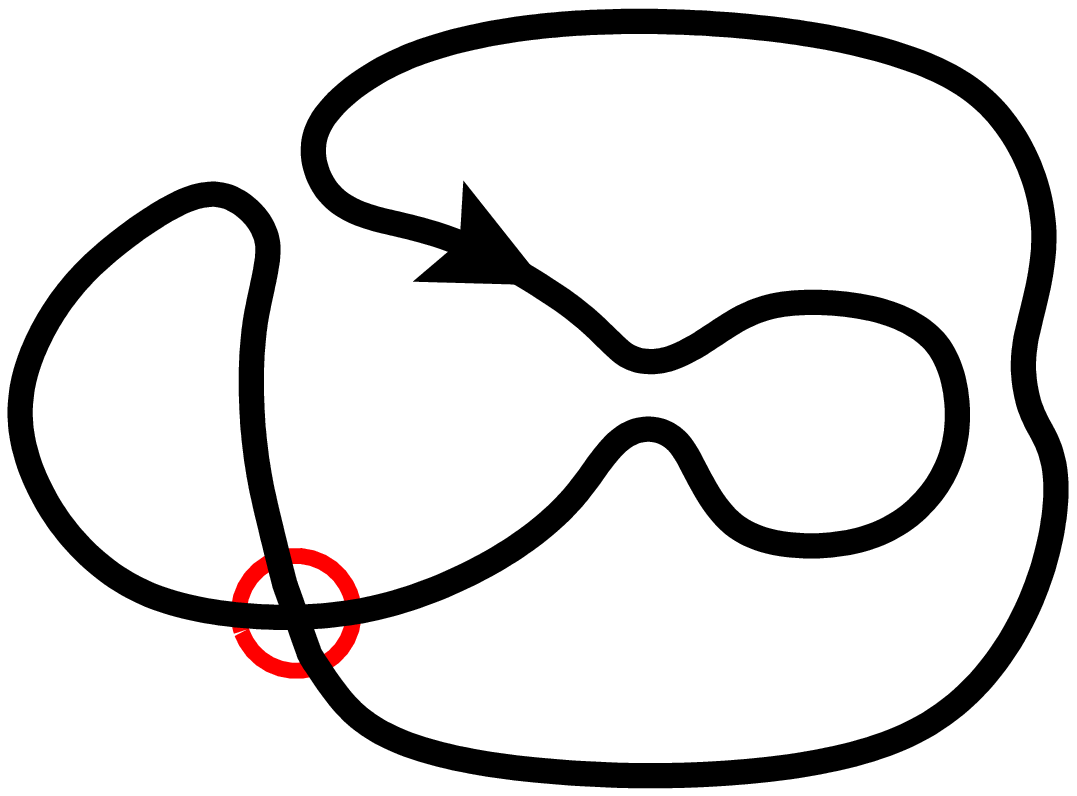} & \laf{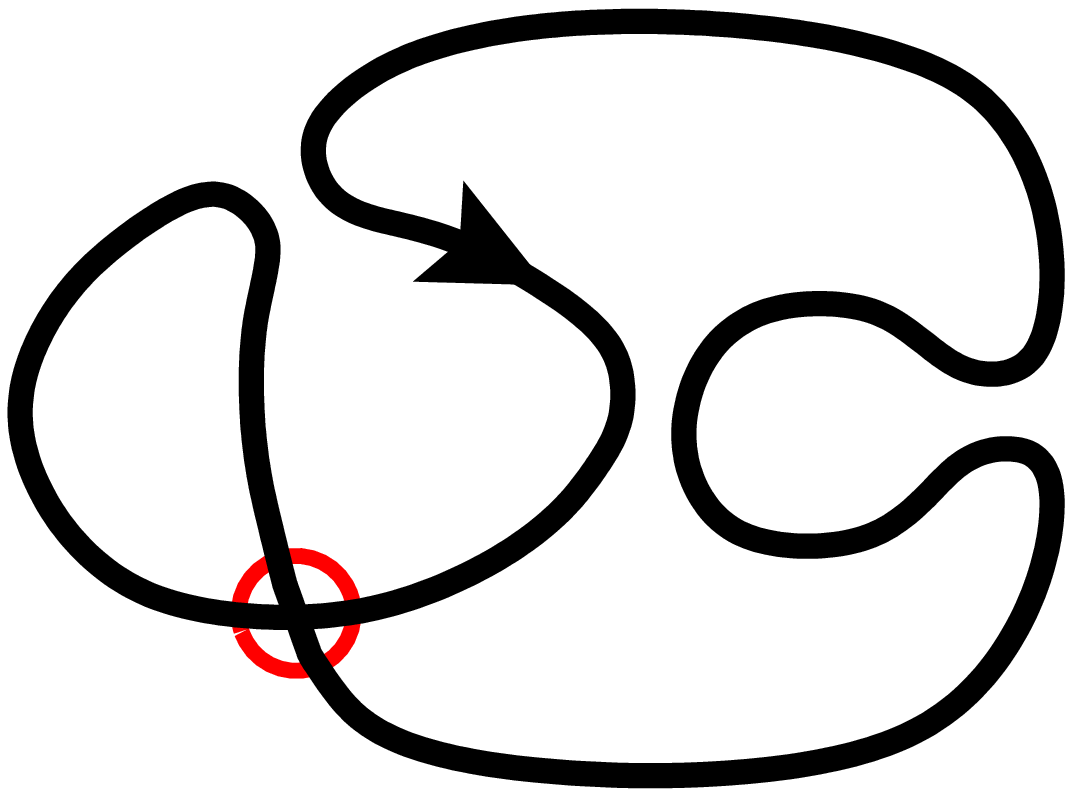} & \laf{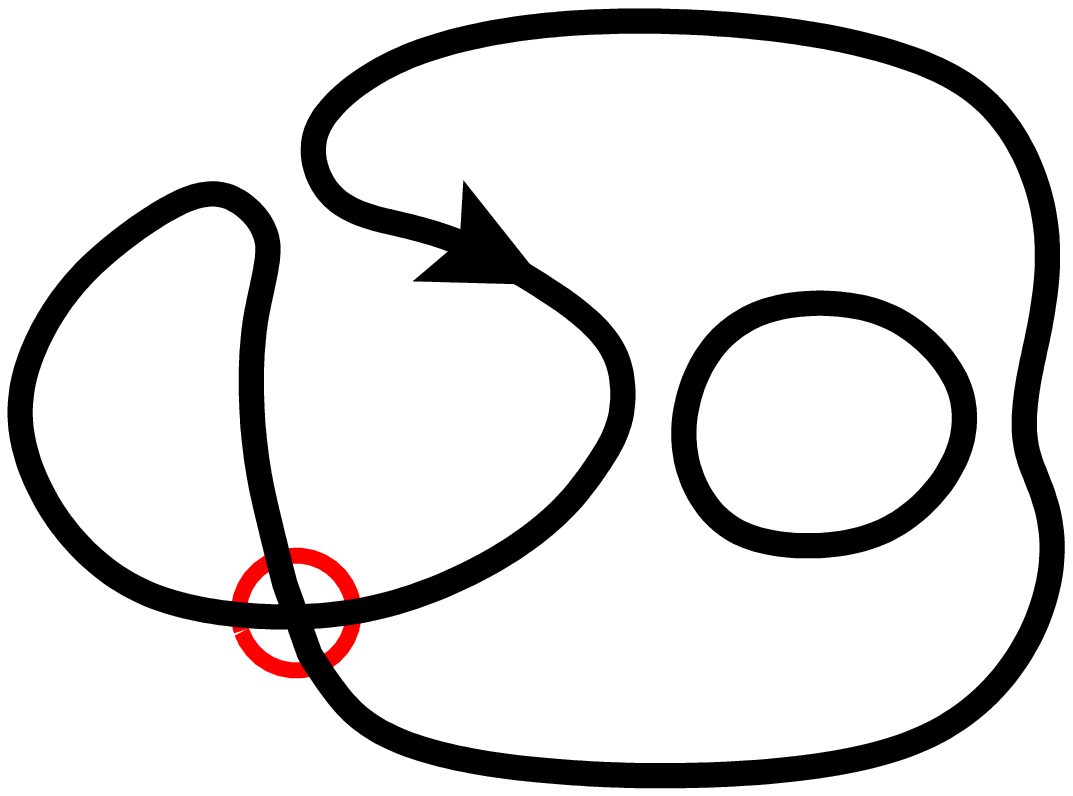}\\ \hline
 (\a,\b,\d) & (3,0,2) & (2,1,1) & (2,1,1) & (1,2,2)\makebox(0,15){}
\\ \hline\hline
& \laf{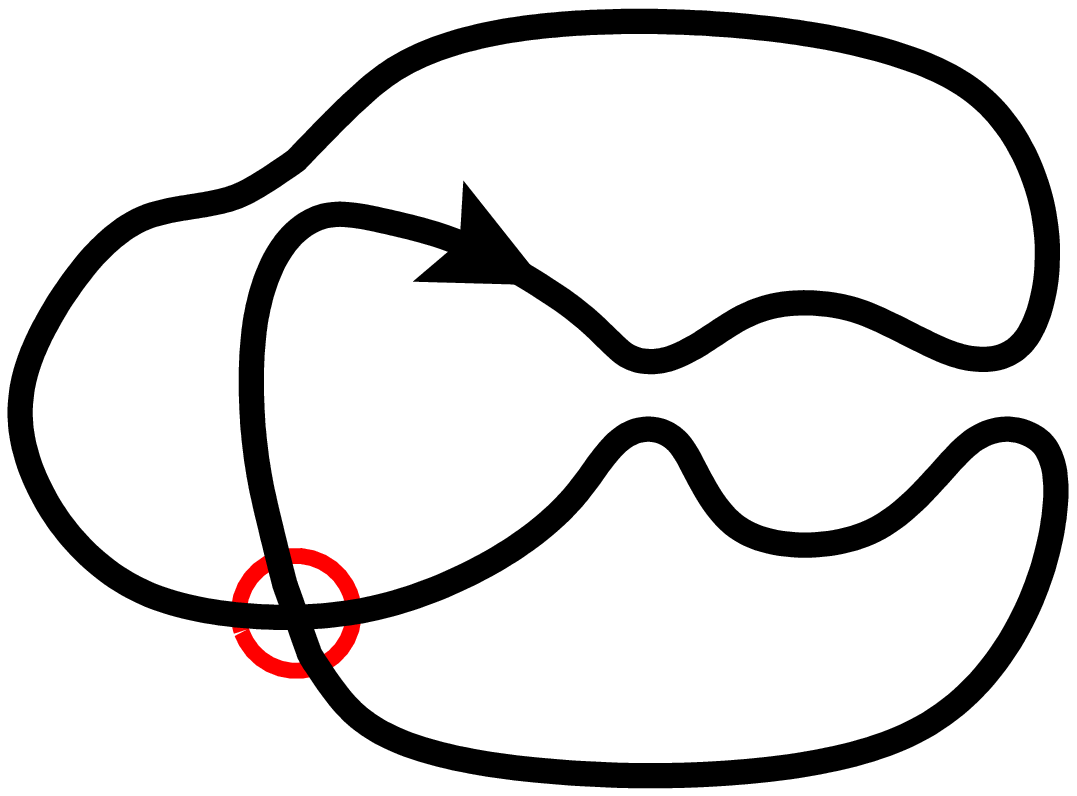} & \laf{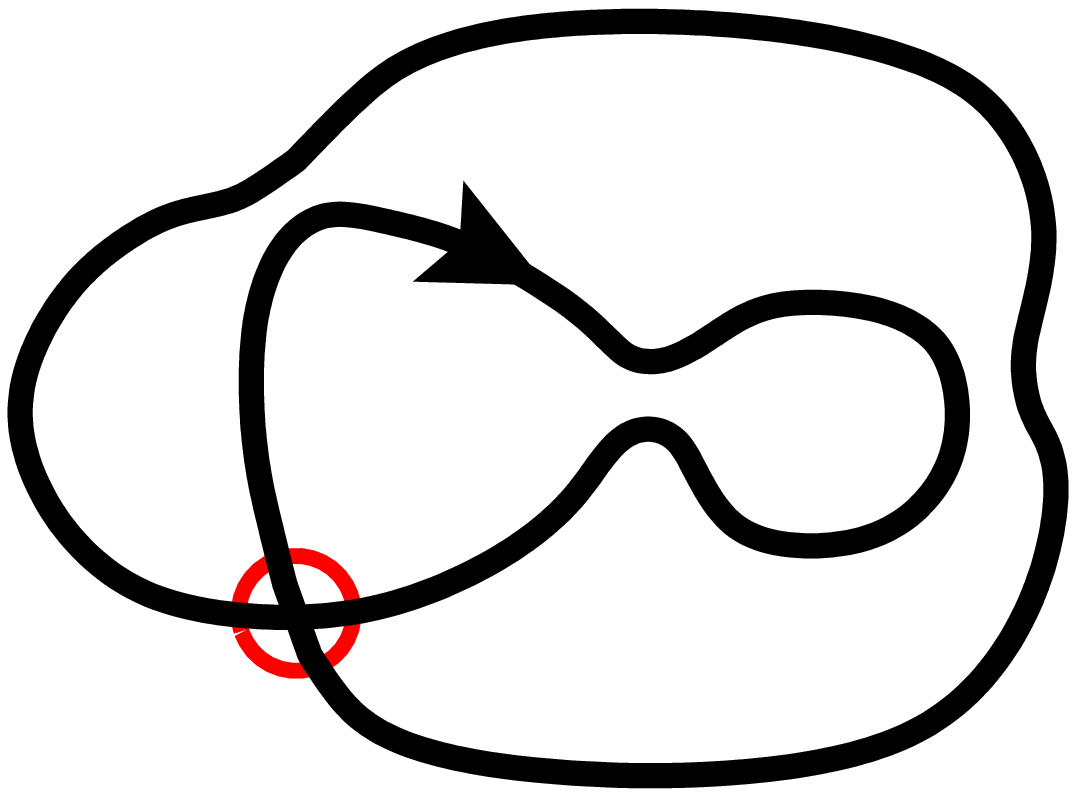} & \laf{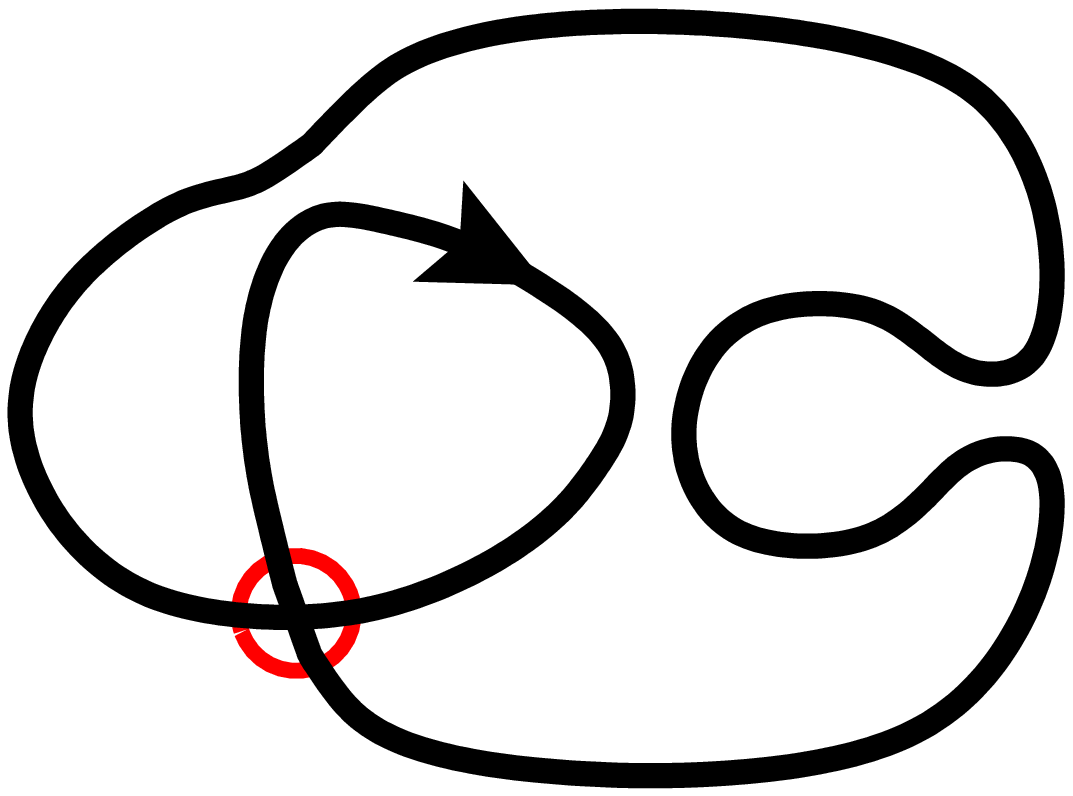} & \laf{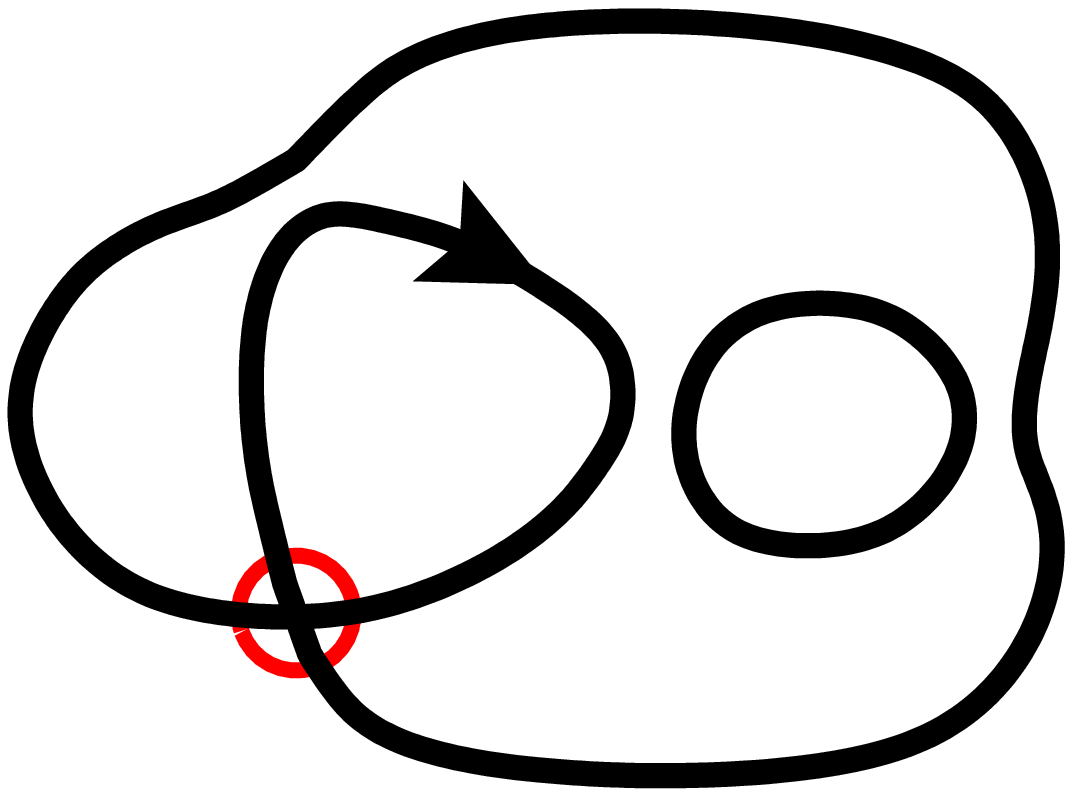}\\ \cline{2-5}
& (2,1,1) & (1,2,1) & (1,2,1) & (0,3,2)\makebox(0,15){}
\end{array}
$$

We have
$$\begin{array}{rcl}
\kb{L} &=& A^3d + 3A^2B + 2AB^2 + AB^2d+ B^3d\ ; \vspace{8pt}\\
J_L(t) &=& t^{-2}-t^{-1}-t^{-1/2}+1+t^{1/2}\ .
\end{array}$$}
\end{Example}

\bigskip
\section{Ribbon graphs and the Bollob\'as-Riordan polynomial}\label{s:rg}

Ribbon graphs are the objects of Topological Graph Theory. There are several books on this subject and its applications \cite{GT,LZ,MT}. Our ribbon graphs are nothing else than
the band decompositions from \cite[Section 3.2]{GT} with the interior of all 2-bands removed. We modify the definitions of \cite{BR3} to signed ribbon graphs.

\begin{defn}\label{def:rb}{\rm
A (signed) {\it ribbon graph} $G$ is a surface (possibly nonorientable) with boundary represented as the union of two sets of closed topological discs called
{\it vertices} $V(G)$ and {\it edges} $E(G)$, satisfying the following
conditions:
\begin{itemize}
\item[$\bullet$] these vertices and edges intersect by disjoint line segments;
\item[$\bullet$] each such line segment lies on the boundary of precisely one vertex and precisely one edge;
\item[$\bullet$] every edge contains exactly two such line segments,
\end{itemize} 
together with a {\it sign function} $\ve: E(G) \to \{\pm 1\}$.}
\end{defn}
Here are a few examples (if the sign is omitted it is assumed to be $+1$).
$$\risS{-18}{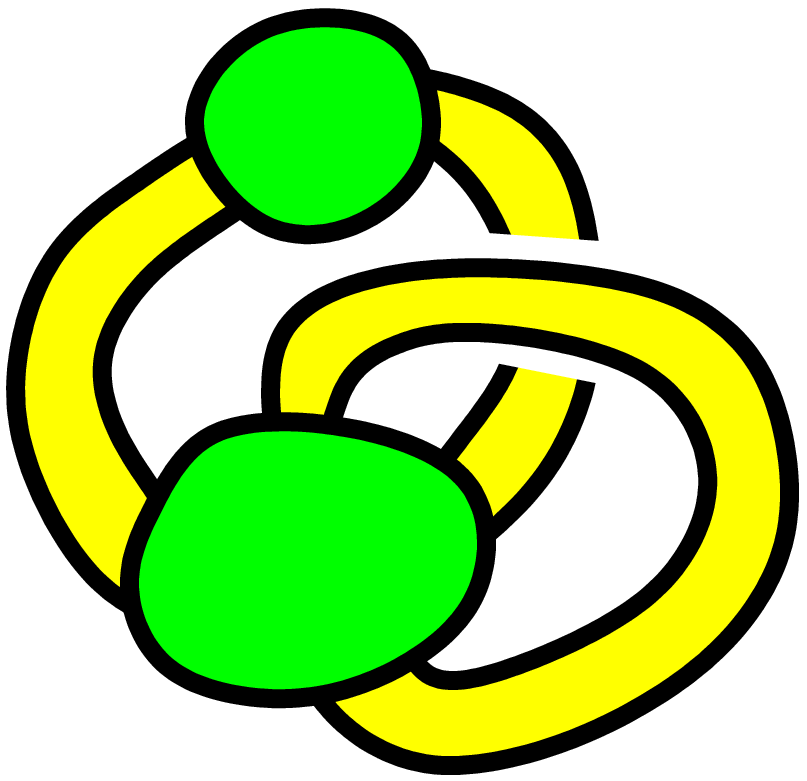}{}{50}{30}{20}\hspace{2cm}
  \risS{-18}{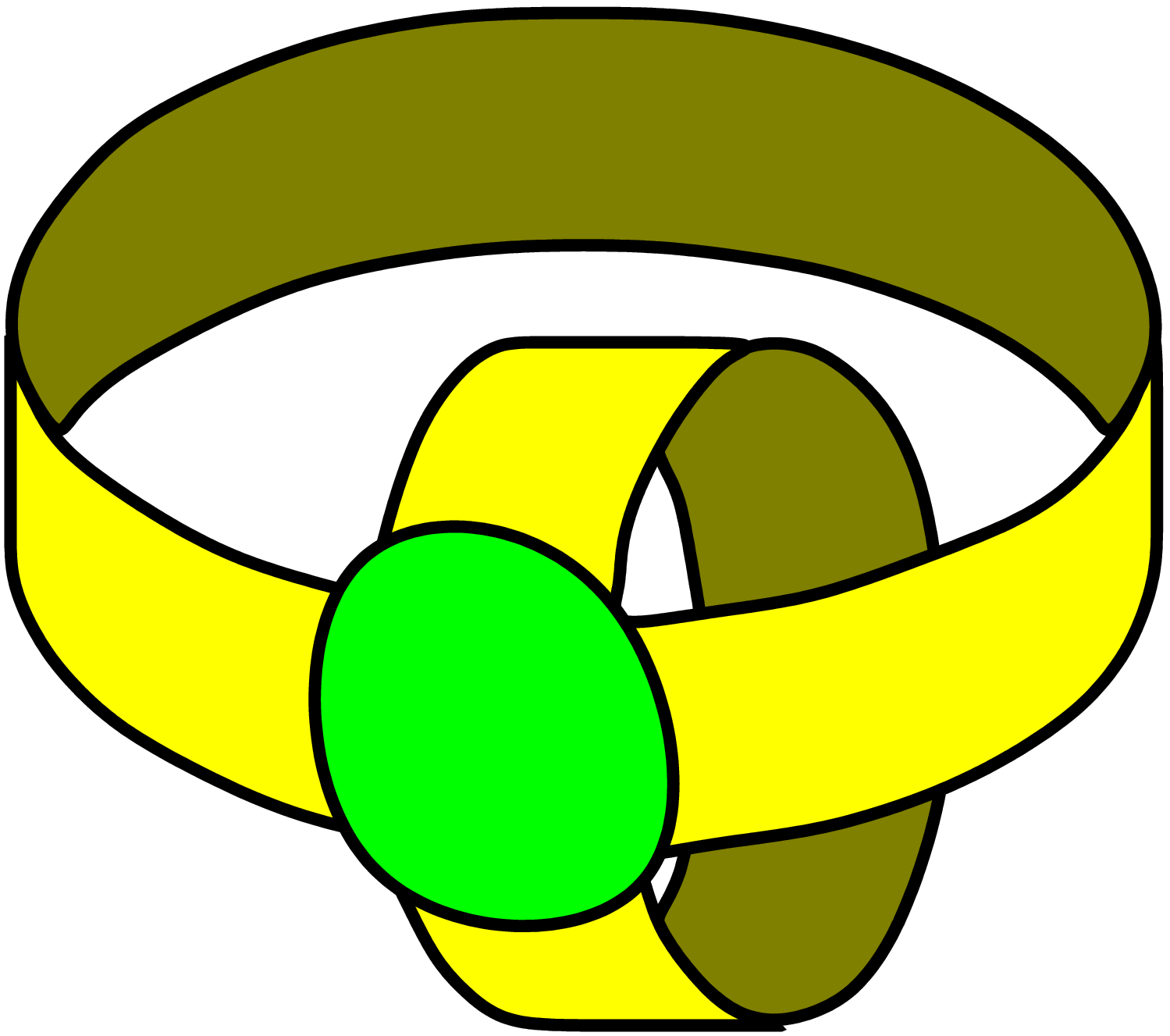}{}{50}{0}{20}\hspace{2cm}
  \risS{-18}{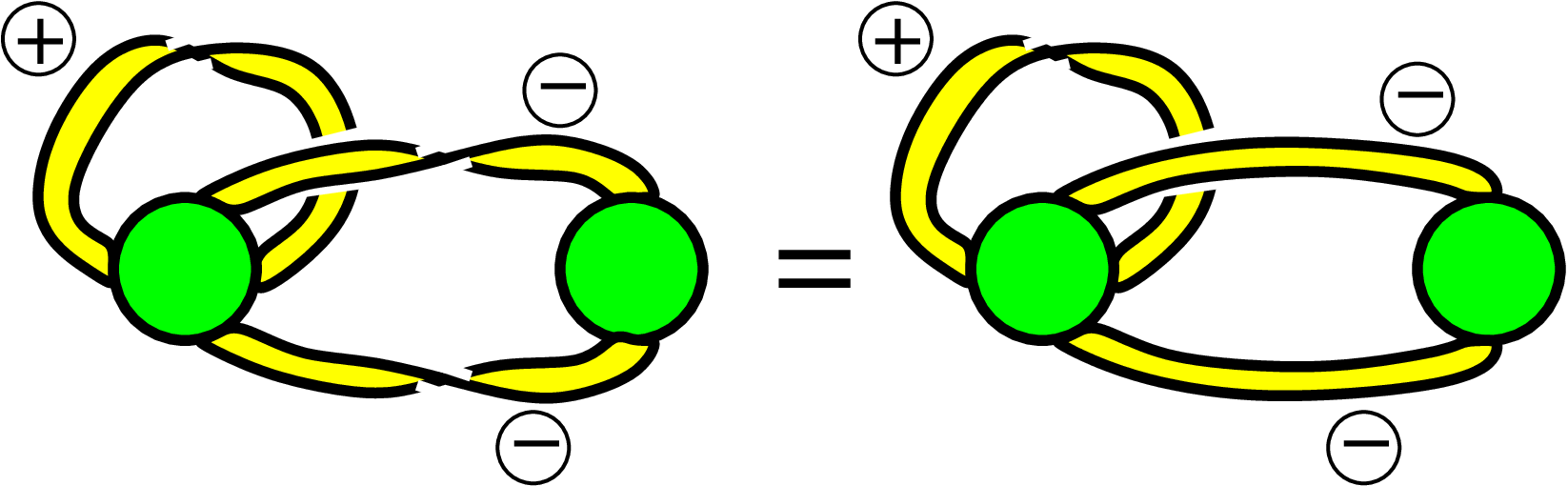}{}{200}{45}{20}
$$
If we put a dot at the center of each vertex-disc and a line in each edge-disc we will get an ordinary graph $\Gamma$, the {\it core graph}, embedded into a surface of $G$. Conversely if we have a graph $\Gamma$ embedded into a surface then it determines a ribbon graph structure on a small neighborhood of $\Gamma$ inside the surface.

To define the Bollob\'as-Riordan polynomial we need to introduce several parameters of  a ribbon graph $G$. Let
\begin{itemize}
\item[$\bullet$] $v(G) := |V(G)|$ denote the number of vertices of $G$;
\item[$\bullet$] $e(G) := |E(G)|$ denote the number of edges of $G$;
\item[$\bullet$] $k(G)$ denote the number of connected components of $G$;
\item[$\bullet$] $r(G):=v(G)-k(G)$ be the {\it rank} of $G$;
\item[$\bullet$] $n(G):=e(G)-r(G)$ be the {\it nullity} of $G$;
\item[$\bullet$] $\bc(G)$ denote the number of connected components of the boundary of the surface of $G$.
\end{itemize} 
A {\it spanning subgraph} of a ribbon graph $G$ is defined as a
subgraph consisting of all the vertices of $G$ and a subset of the edges of $G$.
Let $\cF(G)$ denote the set of spanning subgraphs of $G$.
Clearly, $|\cF(G)| = 2^{e(G)}$. For a signed ribbon graph we need one more parameter of a spanning subgraph. Let $e_{-}(F)$ be the number of negative edges in $F$.
Denote the complement to $F$ in $G$ by $\wo F = G - F$, i.e.
the spanning subgraph of $G$ with exactly those (signed) edges of $G$ that do not belong to $F$.  Finally, let
$$s(F)  = \frac{e_{-}(F)-e_{-}(\wo F)}{2}\ .$$

\begin{defn}\label{def:rb}{\rm
The signed {\it Bollob\'as-Riordan polynomial} $R_G(x,y,z)$ is defined by
$$R_G(x,y,z)\ :=\ \sum_{F \in \cF(G)}
   x^{r(G)-r(F)+s(F)}
   y^{n(F)-s(F)}
   z^{k(F)-\bc(F)+n(F)}\, .
$$}
\end{defn}
In general this is a Laurent polynomial in $x^{1/2}$, $y^{1/2}$, and $z$.

The signed version of the Bollob\'as-Riordan polynomial was introduced in \cite{CP}. If all the edges are positive then it is obtained from the original Bollob\'as-Riordan polynomial
\cite{BR3} by a simple substitution $x+1$ for $x$.
Note that the exponent $k(F)-\bc(F)+n(F)$ of the variable $z$ is equal to
$2k(F)-\chi(\widetilde{F})$, where $\chi(\widetilde{F})$ is the Euler characteristic of the surface $\widetilde{F}$ obtained by gluing a disc to each boundary component of $F$.
For orientable ribbon graphs it is twice the genus of $F$.
In particular, for a planar ribbon graph $G$ (i.e. when the surface~$G$ has genus zero)
the Bollob\'as-Riordan polynomial $R_G$ does not
depend on $z$. In this case it is essentially
equal to the classical Tutte polynomial $T_{\Ga}(x,y)$ of the core graph $\Ga$ of $G$:
$$R_G(x-1,y-1,z) = T_{\Ga}(x,y)$$
if all edges are positive, and if not, to Kauffman's signed Tutte polynomial for signed graphs.  
Similarly, a specialization $z=1$ of the Bollob\'as-Riordan polynomial
of an arbitrary ribbon graph $G$ gives the (signed) Tutte polynomial of the core graph:
$$R_G(x-1,y-1,1) = T_{\Ga}(x,y)\,.$$
So one may think about the Bollob\'as-Riordan polynomial as a generalization of the Tutte polynomial to graphs embedded into a surface.

\begin{Example}\label{ex2} {\rm
Consider the third ribbon graph $G$ from our example above and shown on the left in the table below.
The other columns show eight possible spanning subgraphs $F$ and
the corresponding values of $k(F)$, $r(F)$, $n(F)$, $\bc(F)$, and $s(F)$.
$$
\begin{array}{c||c|c|c|c} \label{br-table}
\laf{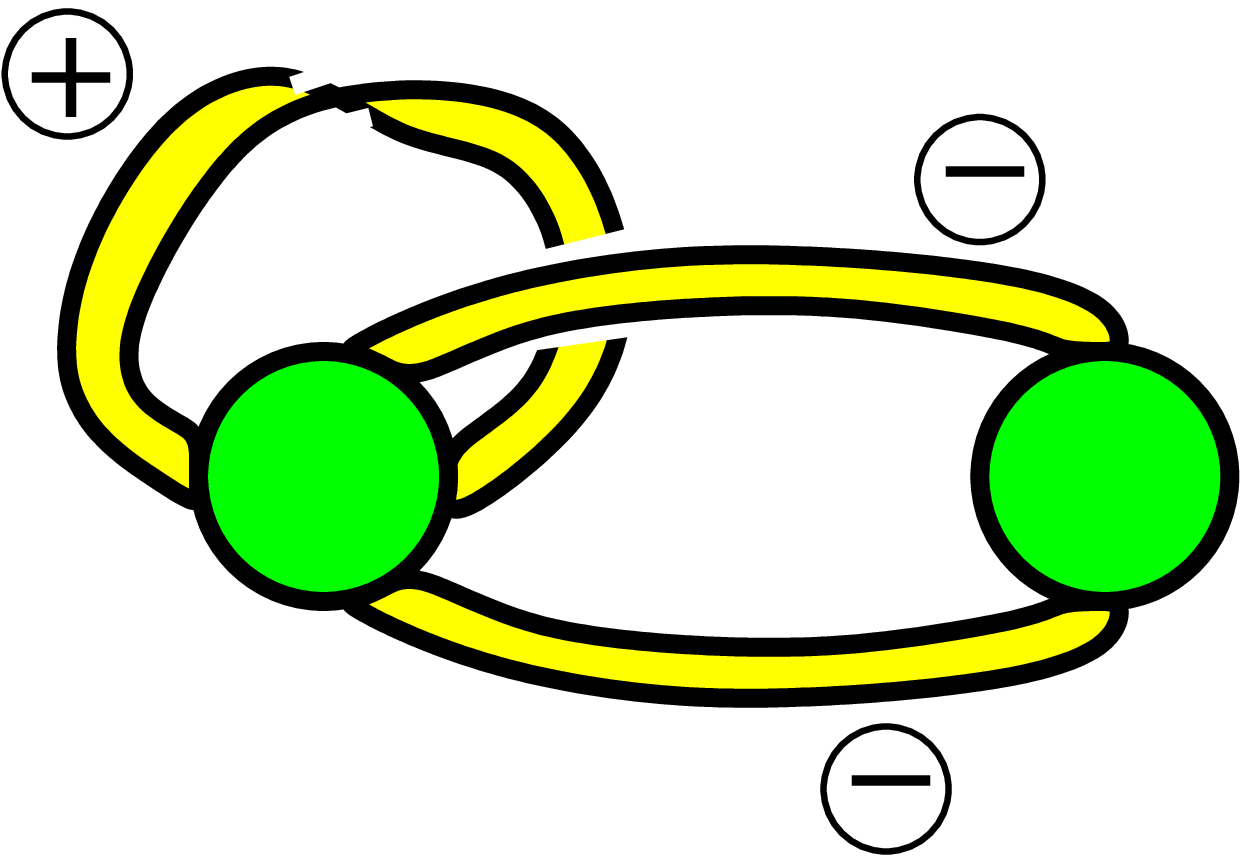} & \quad\risS{0}{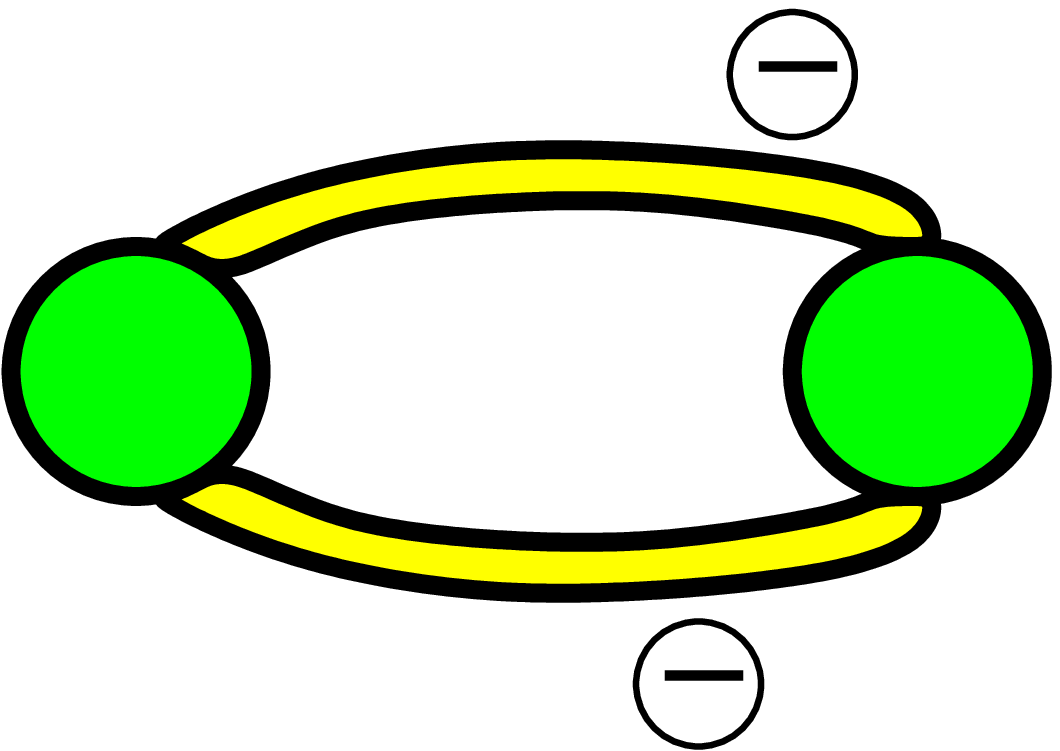}{}{58}{0}{0} & \quad\risS{0}{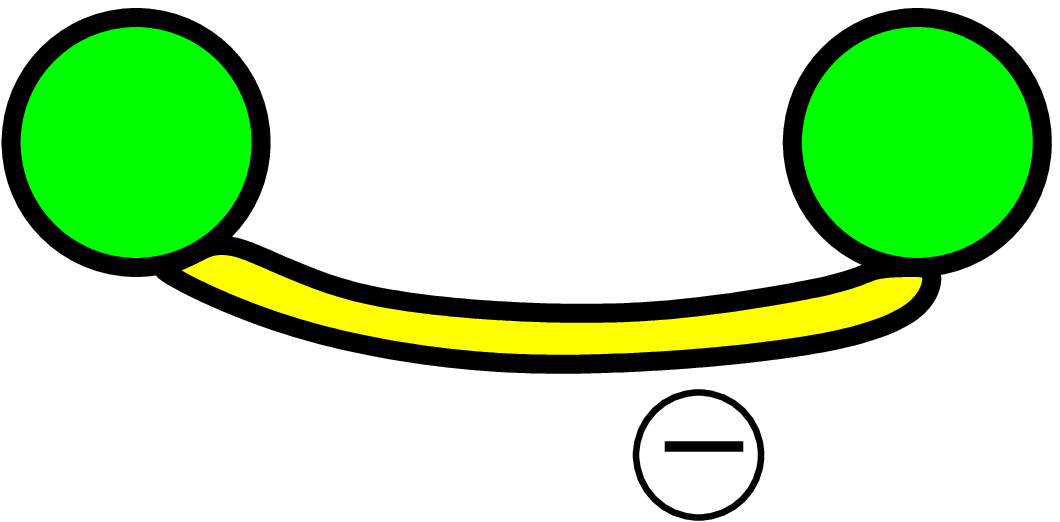}{}{58}{0}{0} & 
\quad\risS{13}{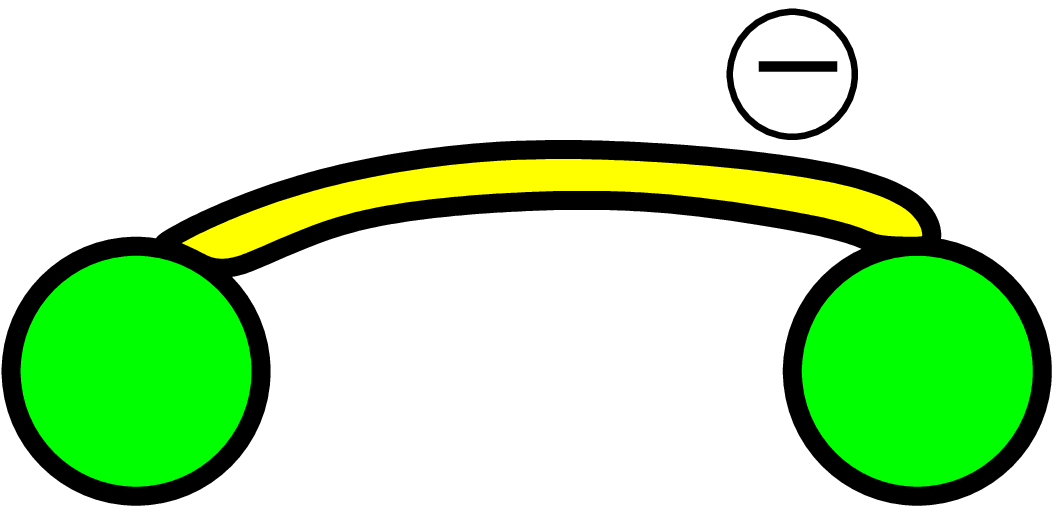}{}{58}{50}{0}  & \quad\risS{13}{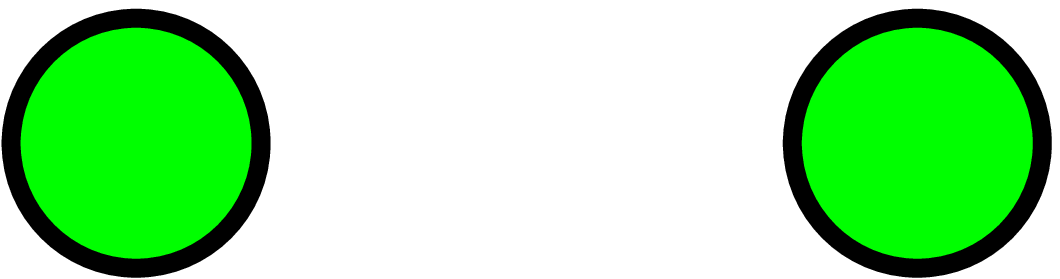}{}{58}{0}{0}\\ \hline
 (k,r,n,\bc,s) & (1,1,1,2,1) & (1,1,0,1,0) & (1,1,0,1,0) & (2,0,0,2,-1)\makebox(0,15){}
\\ \hline\hline
& \laf{rg} & \laf{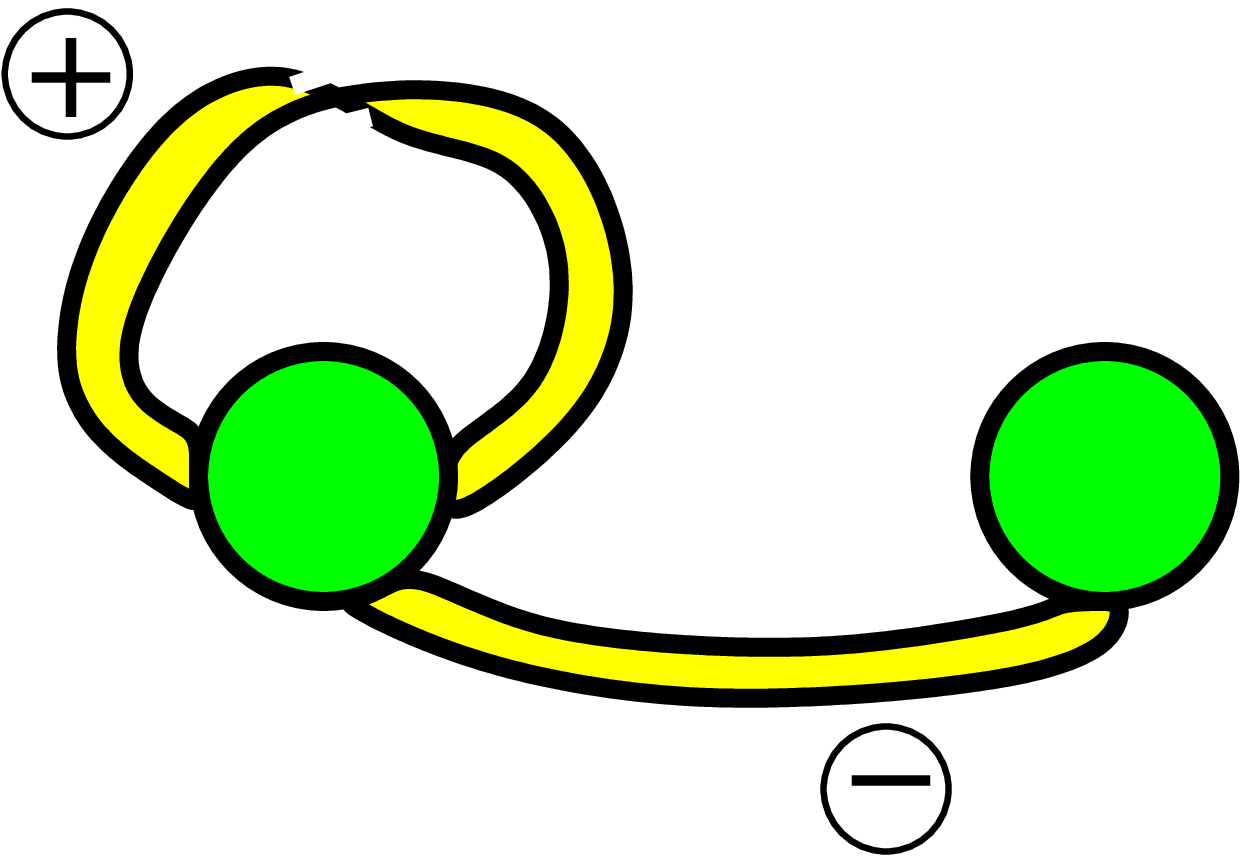} & 
  \risS{13}{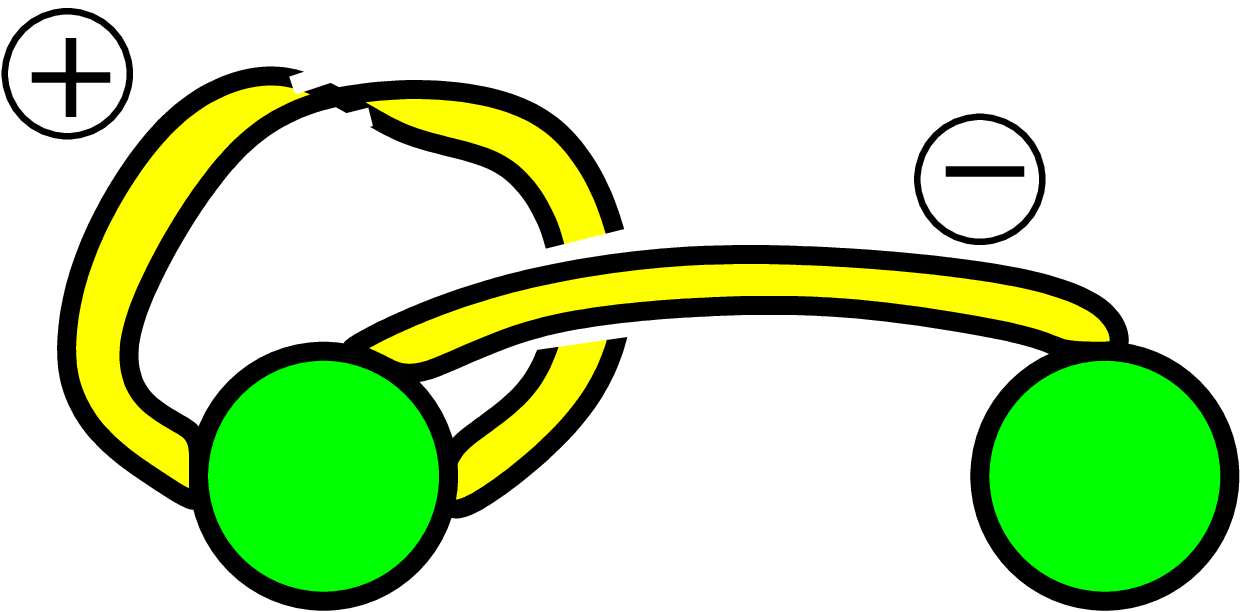}{}{65}{0}{0} & \risS{13}{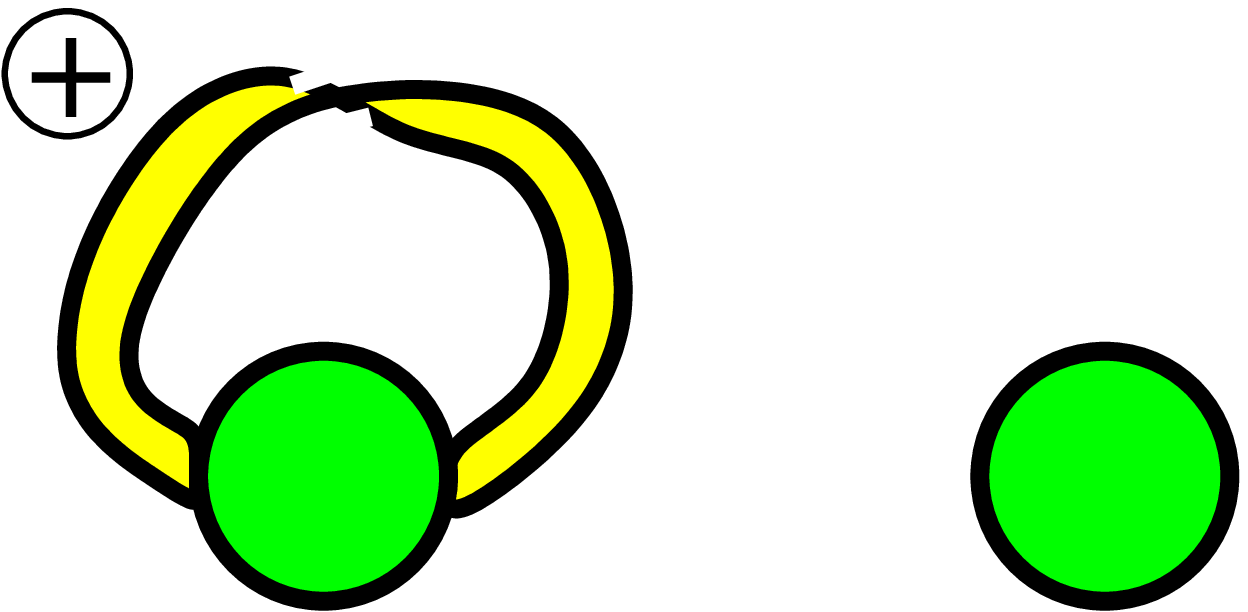}{}{65}{0}{0}\\ \cline{2-5}
& (1,1,2,1,1) & (1,1,1,1,0) & (1,1,1,1,0) & (2,0,1,2,-1)\makebox(0,15){}
\end{array}
$$

We have
$$R_G(x,y,z) = x+2+y+xyz^2+2yz+y^2z\ .$$
}
\end{Example}

\section{Ribbon graphs associated with virtual links}\label{s:vl2rg}
In this section we describe a construction of a ribbon graph starting with a virtual link diagram. Our construction is similar to the classical Seifert algorithm of a construction of the Seifert surface of a link. There are two differences. The first one is that we do not twist the bands in small neighborhoods of the crossings. The second one is that we do not care how our ribbon graph is embedded into the three space.

Suppose our virtual link diagram $L$ is oriented. Then there is a state where all the splittings preserve orientation. Following a suggestion of N.~Stoltzfus, we will call it the {\it Seifert state}, because its state circles are the Seifert circles of the link diagram. Also we will call all splittings in the Seifert state 
{\it Seifert splittings}.

The Seifert circles will be the boundary circles of the vertex-discs of the future ribbon graph. So we are going to glue in a disc to each Seifert circle. Before that, though, let us describe the edges of the ribbon graph. When we are doing a Seifert splitting in a vicinity of a crossing we place a small planar band connecting two branches of the splitting. These bands will be the edge-discs of our ribbon graph. If the Seifert splitting was an $A$-splitting we assign $+1$ to the corresponding edge-band, if it was
a $B$-splitting then we assign $-1$. It is easy to see that this sign is equal to the local writhe of the crossing. So we get a sign function. Because of the presence of virtual crossings our Seifert circles may be twisted, i.e. they are actually immersed into the plane with double points at virtual crossings. In the next step of the construction we untwist all Seifert circles to resolve the double points. This may result in some twisting on the edge-bands. After that we pull all Seifert circles apart, which could lead to additional twisting of our edges. In the last step, we glue the vertex-discs into the circles.
The signed ribbon graph produced is denoted by $G_L$.
The next example illustrates this procedure.

\begin{Example}\label{ex:vl2rg} {\rm
$$\risS{-20}{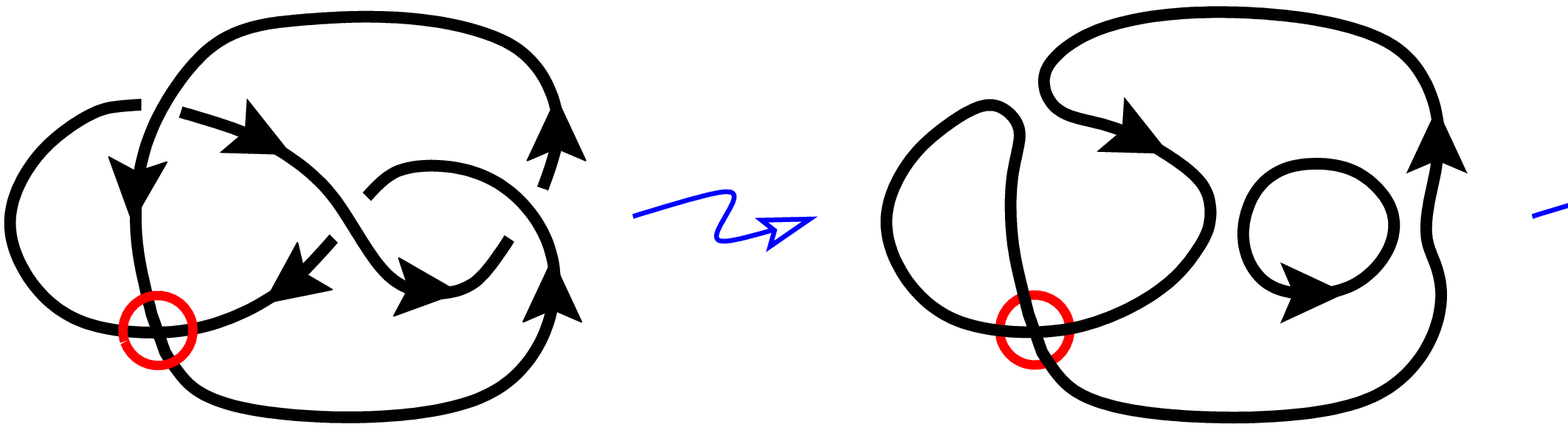}{
       \put(5,-10){\mbox{\scriptsize Diagram}}
       \put(42,33){\mbox{$L$}}
       \put(68,-10){\mbox{\scriptsize Seifert state}}
       \put(125,-10){\mbox{\scriptsize Attaching bands}}
       \put(125,-18){\mbox{\scriptsize to Seifert circles}}
       \put(200,-10){\mbox{\scriptsize Untwisting}}
       \put(200,-18){\mbox{\scriptsize Seifert circles}}
       \put(280,-10){\mbox{\scriptsize Pulling Seifert}}
       \put(280,-18){\mbox{\scriptsize circles apart}}
       \put(370,-10){\mbox{\scriptsize Glue in the}}
       \put(370,-18){\mbox{\scriptsize vertex-discs}}
       \put(385,45){\mbox{$G_L$}}}{420}{25}{50}
$$}
\end{Example}

Another way to explain the same construction is the following. Instead of attaching bands to the Seifert circles we only mark the places on the Seifert circles where the bands have to be attached and memorize the order in which the marks occur according to the orientation on the circles. Then we draw each Seifert circle separately on a plane as a perfect circle oriented counterclockwise. Now attach edge-bands according to the marks;
it is easy to see that this will always result in a half-twist on each band. The sign function is defined as before.

\section{Main Theorem}\label{s:mt}

\begin{thm}\label{th:mt}
Let $L$ be a virtual link diagram, $G_L$ be the corresponding signed ribbon graph, and
$n:=n(G_L)$, $r:=r(G_L)$, $k:=k(G_L)$. Then
$$\kb{L} (A,B,d) = A^n B^r d^{k-1} \,
R_{G_L}\left(\frac{Ad}{B}, \frac{Bd}{A}, \frac{1}{d}\right)\ .
$$
\end{thm}

{\bf Proof.}
Let $L$ be a virtual link diagram, $G_L$ be the corresponding signed ribbon graph, and denote the Seifert state of $G_L$ as $\overline{S}$.
There is a natural bijection between $\cS(L)$, the set of states of $L$, and $\cF(G_L)$, the set of spanning subgraphs of $G_L$.  Namely, given a state $S$, associate to it a spanning subgraph $F_S$ by the following construction.  If a crossing in $S$ is split differently than it is in $\overline{S}$, include its associated edge-band in the spanning subgraph $F_S$.  If a crossing is split the same way for both $S$ and $\overline{S}$, do not include the associated edge-band in $F_S$.  This gives the subgraph $F_S$ associated to $S$.  
(Certainly $|\cS(L)|=|\cF(G_L)|$, because the number of classical crossings of $L$ is equal to $e(G_L)$ by virtue of our construction above.)

\medskip
For example, consider the link $L$ from example \ref{ex1}.  We know from above that $G_L$ is the ribbon graph considered in example \ref{ex2}.  For each state given in the table on page \pageref{kb-table}, we can associate to it a spanning subgraph from the table on page \pageref{br-table}.  Consider the first state $S$ given in the table on page \pageref{kb-table}.  The two rightmost crossings are split differently than they are in the Seifert state of $L$ given above. Thus the spanning subgraph associated to this state $S$ via the correspondence given above is the first subgraph in the table on page \pageref{br-table}.  In fact, each state in the table on page \pageref{kb-table} corresponds correctly to its associated spanning subgraph in the table on page \pageref{br-table}. Check that $\kb{L}$ computed in example \ref{ex1} and
$R_{G_L}(x,y,z)$ computed in example \ref{ex2} satisfy the theorem.

\medskip
Now, given that $F \in \cF(G_L)$ is associated to $S \in \cS(L)$ as described above (for simplicity, we write $F$ instead of $F_S$), consider the term \vspace{8pt}
$$x^{r(G_L)-r(F)+s(F)}
y^{n(F)-s(F)}
z^{k(F)-\bc(F)+n(F)}\ .
\vspace{8pt}
$$ 
Substituting in $x=\frac{Ad}{B}$, $y=\frac{Bd}{A}$, and $z=\frac{1}{d}$ and multiplying by the term $A^nB^rd^{k-1}$ as in the theorem, we have \vspace{8pt}\\
$A^nB^rd^{k-1}(AdB^{-1})^{r-r(F)+s(F)}(BdA^{-1})^{n(F)-s(F)}d^{-k(F)+\bc(F)-n(F)}$ \vspace{2pt}

$$\begin{array}{rcl}
&=& A^{n+r-r(F)-n(F)+2s(F)}B^{r(F)+n(F)-2s(F)}d^{k+r-k(F)-r(F)+\bc(F)-1}\ . \\
\end{array} $$ 
\\
Since $r(G):=v(G)-k(G)$ and $n(G):=e(G)-r(G)$ for any ribbon graph $G$, we can rewrite our term as \vspace{8pt}
$$A^{e(G_L)-e(F)+2s(F)}B^{e(F)-2s(F)}d^{v(G_L)-v(F)+\bc(F)-1}\ .
\vspace{8pt}
$$
And since $v(G_L)=v(F)$ by the definition of a spanning subgraph, we have \vspace{8pt}
\begin{equation}
A^{e(G_L)-e(F)+2s(F)}B^{e(F)-2s(F)}d^{\bc(F)-1}\ .
\vspace{8pt}
\end{equation}
It suffices to show that this is equal to the Kauffman bracket term $A^{\a(S)}\,B^{\b(S)}\,d^{\,\d(S)-1}$, since our bijection described above will then imply the theorem.  
We first show that $e(F)-2s(F)=\b(S)$ by a counting argument.  Using the definition of $s(F)$, we get
\begin{equation}
e(F)-2s(F)=e(F)-e_{-}(F)+e_{-}(\overline{F})\ .
\end{equation} 

Consider the crossings of $L$ and how they are split in $\overline{S}$.  Let $m$ denote the number of crossings which are $B$-splittings in $\overline{S}$.  Let $b$ denote the number of crossings which are $B$-splittings in $\overline{S}$ but are $A$-splittings in $S$. Let $a$ denote the number of crossings which are $A$-splittings in $\overline{S}$ but are $B$-splittings in $S$.

Now, since $e(F)$ is the number of edges included in $F$, $e(F)$ equals the number of crossings of $L$ which are split differently between $S$ and $\overline{S}$.  That is, $e(F)=a+b$.  
Recall that $e_{-}(F)$ denotes the number of edges in $F$ with sign $-1$.  And since any such edge corresponds to a $B$-splitting in $\overline{S}$, it is clear that $e_{-}(F)=b$.  

Since $\overline{F}$ is the complement of $F$, $e_{-}(\overline{F})$ denotes the number of crossings which are $B$-splittings in $\overline{S}$ and also in $S$.  So, we deduce that $e_{-}(\overline{F})=m-b$.  
Finally, we consider $\b(S)$, the number of crossings which are $B$-splittings in $S$.  So clearly, $\b(S)=a+e_{-}(\overline{F})=a+(m-b)$.  

Thus, we have that 
$$\begin{array}{rcl}
e(F)-e_{-}(F)+e_{-}(\overline{F}) &=& (a+b)-b+(m-b)\\
&=& a+(m-b)\\
&=& \b(S)\ . 
\end{array}$$
This with $(2)$ gives the desired result.  And the fact that $e(G_L)-e(F)+2s(F)=\a(S)$ is immediate, since we just showed that $e(F)-2s(F)=\b(S)$, and certainly $e(G_L)=\a(S)+\b(S)$. 

To finish the proof it remains to show that $\bc(F)=\d(S)$. For that let us trace simultaneously a circle of the state $S$ and a boundary component of $F$. Suppose we are passing a place near a crossing. If this crossing is split in the same way as in $\wo{S}$ then we continue to go along the circle which is locally the same as in $\wo{S}$. On the whole ribbon graph $G_L$ this means that we are passing a place where the corresponding edge-band is supposed to be attached. Or, in other words, we skip the edge and do not include it into the corresponding spanning subgraph. This is precisely how we obtained $F$. Now if the passed crossing in $S$ is split differently compared with $\wo{S}$, then we have to switch to another strand of $\wo{S}$. This means that we should turn on the corresponding edge-band in $G_L$, i.e. we should include this edge-band in the corresponding subgraph. Again this is precisely what we did with $F$. The next picture with the first state of example \ref{ex1} (the table on page \pageref{kb-table}) illustrates this.
$$\risS{-20}{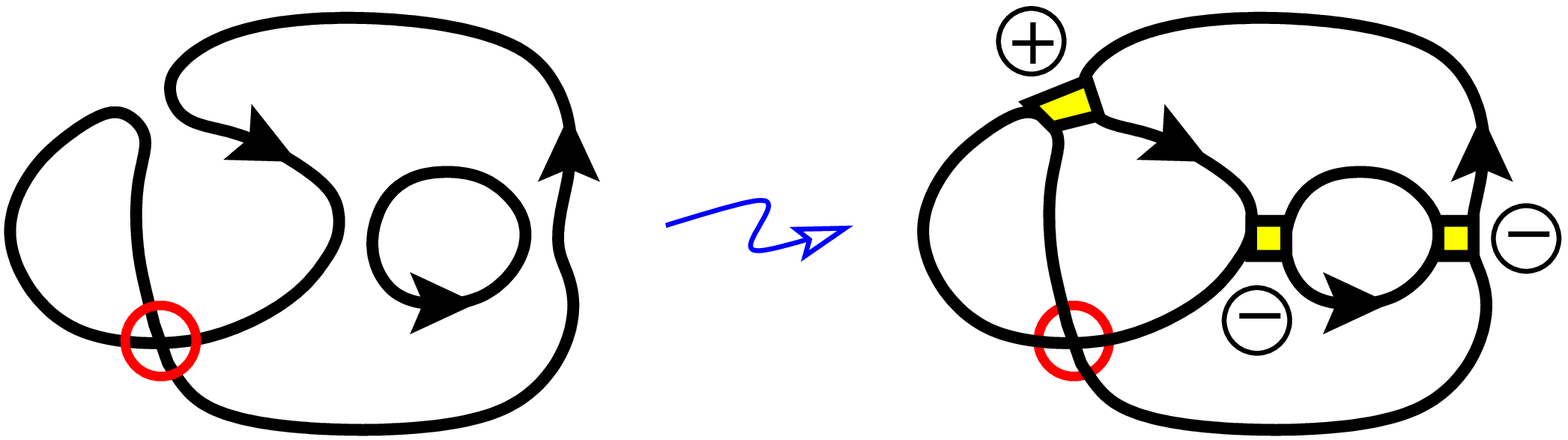}{
       \put(-5,-10){\mbox{\rm\scriptsize Seifert state $\wo{S}$}}
       \put(95,-10){\mbox{\scriptsize $G_L$}}}{125}{20}{40}
  \hspace{3cm}
  \risS{-20}{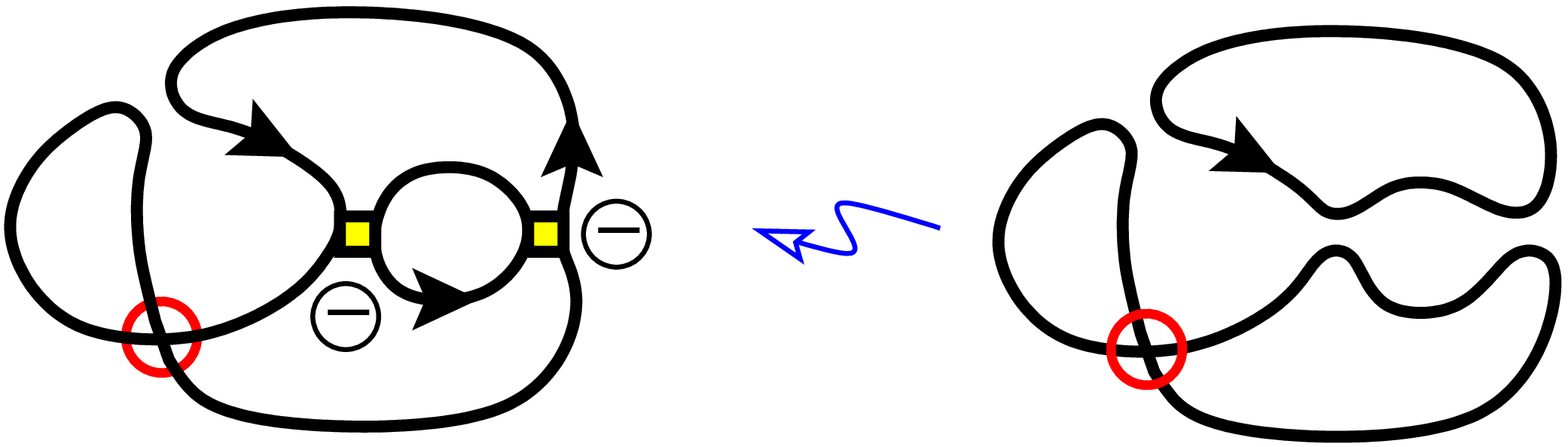}{
       \put(-25,-10){\mbox{\rm\scriptsize Spanning subgraph $F=F_S$}}
       \put(100,-10){\mbox{\scriptsize $S$}}}{125}{20}{35}
$$
Therefore, the tracing of the state circles of $S$ corresponds to the tracing of the boundary components of $F$, i.e. $\bc(F)=\d(S)$.

So, we have shown that $(1)$ is equal to the term of $\kb{L}$ corresponding to the state $S$, and thus theorem \ref{th:mt} is proved.  
\hspace{\fill}$\square$

\vskip1.cm

\parbox[t]{2.5in}{\it \textbf{Sergei~Chmutov}\\
Department of Mathematics\\
The Ohio State University, Mansfield\\
1680 University Drive\\
Mansfield, OH 44906\\
~\texttt{chmutov@math.ohio-state.edu}} \qquad\qquad
\parbox[t]{2.5in}{\it \textbf{Jeremy Voltz}\\
Department of Mathematics,\\
The Ohio State University, \\
231 W. 18th Avenue, \\
Columbus, Ohio 43210\\
~\texttt{voltz.10@osu.edu}}

\end{document}